\documentclass[11pt,oneside,british]{article}
\newcommand{\comment}[1]		{}
\usepackage{amsmath,amssymb,amsthm,rotating,pdflscape,afterpage,scalerel}
\usepackage[sc,osf]{mathpazo}
\usepackage[hyphens]{url}
\AtBeginDocument{
	\DeclareSymbolFont{AMSb}{U}{msb}{m}{n}
	\DeclareSymbolFontAlphabet{\mathbb}{AMSb}
}
\newcommand{\mockalph}[1]{\!}
\usepackage{enumitem}
\usepackage[nouppercase]{scrlayer-scrpage}
\usepackage{avant,mathrsfs,mathtools,faktor,xfrac}
\usepackage[all]{xypic}\xyoption{rotate}
\newdir{ >}{{}*!/-7pt/\dir{>}}
\newdir{i}{{}*!/-7pt/\dir{(}	}
\usepackage{thmtools}
\usepackage[utf8]{inputenc}
\usepackage[T1]{fontenc}
\usepackage[main=english]{babel}
\usepackage{tocloft}
\makeatletter
\renewcommand{\l@figure}{\@dottedtocline{1}{1em}{3.5em}}
\renewcommand{\l@table}{\@dottedtocline{2}{1em}{3.5em}}
\makeatother
\newcommand*{\noaddvspace}{\renewcommand*{\addvspace}[1]{}}
\addtocontents{lof}{\protect\noaddvspace}

\usepackage[hyphens]{url}
\usepackage{breakurl}

\usepackage[
bookmarksopen,
bookmarksdepth=2,
]{hyperref}	
\usepackage[nameinlink]{cleveref} 	

\usepackage{mathabx}

\makeatletter
\let\c@figure\c@table
\makeatother

\numberwithin{table}{section}
\numberwithin{figure}{section}
\newtheorem{abcthm}[table]{Theorem}

\newtheorem{theorem}[table]{Theorem}
\newtheorem{thmx}{Theorem}

\newtheorem{proposition}[table]{Proposition}

\newtheorem{corollary}[table]{Corollary}
\newtheorem{lemma}[table]{Lemma}

\newtheorem{claim}[table]{Claim}

\theoremstyle{definition}

\newtheorem{definition}[table]{Definition}

\newtheorem{notation}[table]{Notation}

\newtheorem{conjecture}[table]{Conjecture}

\theoremstyle{remark}
\newtheorem{fact}[table]{Fact}

\newtheorem{example}[table]{Example}
\newtheorem{exercise}[table]{Exercise}

\newtheorem{histrmks}[table]{Historical remarks}
\newtheorem{remark}[table]{Remark}
\newtheorem{remarks}[table]{Remarks}

\theoremstyle{plain}
\newtheorem*{thm*}{Theorem}
\newtheorem*{theorem*}{Theorem}
\newtheorem*{prop*}{Proposition}
\newtheorem*{proposition*}{Proposition}
\newtheorem*{lemma*}{Lemma}
\newtheorem*{corollary*}{Corollary}
\newtheorem*{cor*}{Corollary}

\theoremstyle{definition}
\newtheorem*{definition*}{Definition}
\newtheorem*{defn*}{Definition}
\newtheorem*{QQ*}{Question}
\newtheorem*{obs*}{Observation}
\newtheorem*{notation*}{Notation}

\theoremstyle{remark}
\newtheorem*{rmk*}{Remark}
\newtheorem*{remark*}{Remark}

\newtheorem*{examples*}{Examples}
\newtheorem*{example*}{Example}
\newtheorem*{EG*}{Example}
\newtheorem*{EGs*}{Examples}
\newtheorem*{fact*}{Fact}
\newtheorem*{prob*}{Problem}

\usepackage{etoolbox}
\usepackage[usenames,dvipsnames,svgnames,table]{xcolor}
\newcommand		{\defd}[1]	{\textcolor{RoyalBlue}{\textbf{\textit{#1}}}}
\newcommand		{\defm}[1]	{\textcolor{RoyalBlue}{#1}}

\tracingpatches
\makeatletter
\patchcmd{\@setref}{\bfseries ??}{\bfseries\color{red} FIX ME!}{}{}
\patchcmd{\@setcite}{\bfseries ?}{\bfseries\color{red} FIX ME!}{}{}
\patchcmd{\@setcref}         {??}{\color{red} FIX ME!}{}{}
\patchcmd{\@setcref}         {??}{\color{red} FIX ME!}{}{}
\patchcmd{\@setcrefrange}    {??}{\color{red} FIX ME!}{}{}
\patchcmd{\@setcrefrange}    {??}{\color{red} FIX ME!}{}{}
\patchcmd{\@setcrefrange}    {??}{\color{red} FIX ME!}{}{}
\patchcmd{\@setcrefrange}    {??}{\color{red} FIX ME!}{}{}
\patchcmd{\@setcrefrange}    {??}{\color{red} FIX ME!}{}{}
\patchcmd{\@setcrefrange}    {??}{\color{red} FIX ME!}{}{}
\patchcmd{\@setnamecref}     {??}{\color{red} FIX ME!}{}{}
\patchcmd{\@setnamecref}     {??}{\color{red} FIX ME!}{}{}
\patchcmd{\@setcpageref}     {??}{\color{red} FIX ME!}{}{}
\patchcmd{\@setcpageref}     {??}{\color{red} FIX ME!}{}{}
\patchcmd{\@setcpagerefrange}{??}{\color{red} FIX ME!}{}{}
\patchcmd{\@setcpagerefrange}{??}{\color{red} FIX ME!}{}{}
\patchcmd{\@setcpagerefrange}{??}{\color{red} FIX ME!}{}{}
\patchcmd{\@setcpagerefrange}{??}{\color{red} FIX ME!}{}{}
\patchcmd{\@setcpagerefrange}{??}{\color{red} FIX ME!}{}{}
\patchcmd{\@cref}            {??}{\color{red} FIX ME!}{}{}
\def\blx@citation@entry#1#2{%
	\blx@bibreq{#1}%
	\ifinlist{#1}{\blx@cites}
	{}
	{\listgadd{\blx@cites}{#1}%
		\blx@auxwrite\@mainaux{}{\string\abx@aux@cite{#1}}}%
	\ifinlistcs{#1}{blx@segm@\the\c@refsection @\the\c@refsegment}
	{}
	{\listcsgadd{blx@segm@\the\c@refsection @\the\c@refsegment}{#1}}%
	\blx@ifdata{#1}%
	{}%
	{\ifcsdef{blx@miss@\the\c@refsection}%
		{\ifinlistcs{#1}{blx@miss@\the\c@refsection}%
			{{\bfseries\color{red} cite:} }%
			{\blx@logreq@active{#2{#1}}}}%
		{\blx@logreq@active{#2{#1}}}}}
\def\blx@citeadd#1{%
	\ifcsdef{blx@keyalias@\the\c@refsection @#1}
	{\edef\blx@realkey{\csuse{blx@keyalias@\the\c@refsection @#1}}}
	{\def\blx@realkey{#1}}%
	\expandafter\blx@citation\expandafter{\blx@realkey}\blx@msg@cundefon
	\expandafter\blx@ifdata\expandafter{\blx@realkey}
	{\advance\blx@tempcnta\@ne
		\listeadd\blx@tempa{\blx@realkey}}
	{\ifnum\blx@tempcntb>\z@\multicitedelim\fi
		\expandafter\abx@missing\expandafter{\blx@realkey}%
		\advance\blx@tempcntb\@ne}}

\DeclarePairedDelimiterX{\pmodx}[1]{(}{)}{{\operator@font mod}\mkern6mu#1}
\renewcommand{\pmod}{%
	\allowbreak
	\if@display\mkern18mu\else\mkern8mu\fi
	\pmodx
}
\makeatother

\makeatletter
\newcommand{\oset}[3][0ex]{%
	\raisebox{.175ex}{$%
		\mathrel{\mathop{#3}\limits^{
				\vbox to#1{\kern-2\ex@
					\hbox{$\scriptstyle#2$}\vss}}}
		$}%
}
\makeatother

\newcommand{\myred}{BrickRed}
\hypersetup{
	colorlinks			= true, 	
	urlcolor			= \myred, 
	linkcolor			= \myred, 	
	citecolor			= PineGreen 	
}

\usepackage{xspace}

\usepackage[left=2.54cm,right=2.54cm,top=2.54cm,bottom=2.54cm]{geometry}
\usepackage{tikz}
\usetikzlibrary{matrix,fit,backgrounds,arrows.meta,calc}
\tikzset{>={Stealth[length=3.6pt,inset=2.35pt,width=4.75pt,round]}}

\makeatletter
\newbox\xrat@below
\newbox\xrat@above
\newcommand{\xrightarrowtail}[2][]{%
	\setbox\xrat@below=\hbox{\ensuremath{\scriptstyle #1}}%
	\setbox\xrat@above=\hbox{\ensuremath{\scriptstyle #2}}%
	\pgfmathsetlengthmacro{\xrat@len}{max(\wd\xrat@below,\wd\xrat@above)+.6em}%
	\mathrel{\tikz [>->,baseline=-.55ex]
		\draw (0,0) -- node[below=-2pt] {\box\xrat@below}
		node[above=-2pt] {\box\xrat@above}
		(\xrat@len,0) ;}}
\newbox\xrat@below
\newbox\xrat@above
\newcommand{\xtwoheadrightarrow}[2][]{%
	\setbox\xrat@below=\hbox{\ensuremath{\scriptstyle #1}}%
	\setbox\xrat@above=\hbox{\ensuremath{\scriptstyle #2}}%
	\pgfmathsetlengthmacro{\xrat@len}{max(\wd\xrat@below,\wd\xrat@above)+.6em}%
	\mathrel{\tikz [-{>[sep=-.5pt]>},baseline=-.55ex]
		\draw (0,0) -- node[below=-2pt] {\box\xrat@below}
		node[above=-2pt] {\box\xrat@above}
		(\xrat@len,0) ;}}
\makeatother
\newcommand{\xmono}{\xrightarrowtail}
\newcommand{\mono}{\xmono{\phantom{\ \, }}}
\newcommand{\xepi}{\xtwoheadrightarrow}
\newcommand{\epi}{\xepi{\phantom{\ \, }}}

\makeatletter
\newcommand{\presectionskip}{-1.5\baselineskip}
\newcommand{\postsectionskip}{0.3\baselineskip}
\usepackage{titlesec}
\renewcommand{\section}{\@startsection
	{chapter}{0}{0mm}
	{\presectionskip}
	{\postsectionskip}
	{\sffamily\huge}}
\renewcommand{\section}{\@startsection
	{section}{1}{0mm}
	{\presectionskip}
	{\postsectionskip}
	{\sffamily\LARGE}}
\renewcommand{\subsection}{\@startsection
	{subsection}{2}{0mm}
	{\presectionskip}
	{\postsectionskip}
	{\sffamily\Large}}
\renewcommand{\subsubsection}{\@startsection
	{subsubsection}{3}{0mm}
	{\presectionskip}
	{\postsectionskip}
	{\sffamily\normalsize}}
\renewcommand{\@seccntformat}[1]{\csname the#1\endcsname.\quad}
\newcommand\HUGE{\@setfontsize\Huge{30}{47}} 

\titleformat{\chapter}[display]
{\sffamily\Large}
{Chapter {\HUGE\normalfont\thechapter}}    
{1em}
{\huge}
\titleformat{\section}
{\sffamily\Large}
{\normalfont{\@setfontsize\Large{18}{47}\thesection.}}    
{1em}
{}
\titleformat{\subsection}
{\sffamily\large}
{\normalfont{\Large\thesubsection}.}    
{1em}
{}
\titleformat{\subsubsection}
{\sffamily\normalsize}
{\normalfont{\large\thesubsubsection}.}    
{1em}
{}
\makeatother

\makeatletter
\def\smallunderbrace#1{\mathop{\vtop{\m@th\ialign{##\crcr
				$\hfil\displaystyle{#1}\hfil$\crcr
				\noalign{\kern3\p@\nointerlineskip}%
				\tiny\upbracefill\crcr\noalign{\kern3\p@}}}}\limits}
\makeatother

\renewcommand{\SS}{\textsection}
\newcommand{\bthm}{\begin{theorem}}
	\newcommand{\ethm}{\end{theorem}}
\newcommand{\bprop}{\begin{proposition}}
	\newcommand{\eprop}{\end{proposition}}
\newcommand{\bcor}{\begin{corollary}}
	\newcommand{\ecor}{\end{corollary}}
\newcommand{\bconj}{\begin{conjecture}}
	\newcommand{\econj}{\end{conjecture}}
\newcommand{\blem}{\begin{lemma}}
	\newcommand{\elem}{\end{lemma}}
\newcommand{\bclm}{\begin{claim}}
	\newcommand{\eclm}{\end{claim}}
\newcommand{\bpf}{\begin{proof}}
	\newcommand{\epf}{\end{proof}}
\newcommand{\bdetails}{\begin{details}}
	\newcommand{\edetails}{\end{details}}
\newcommand{\bdefi}{\begin{definition}}
	\newcommand{\edefi}{\end{definition}}
\newcommand{\bdefn}{\begin{definition}}
	\newcommand{\edefn}{\end{definition}}
\newcommand{\bex}{\begin{example}}
	\newcommand{\eex}{\end{example}}
\newcommand{\bprob}{\begin{problem}}
	\newcommand{\eprob}{\end{problem}}
\newcommand{\bexer}{\begin{exercise}}
	\newcommand{\eexer}{\end{exercise}}
\newcommand{\bexers}{\begin{exercises}}
	\newcommand{\eexers}{\end{exercises}}
\newcommand{\brmk}{\begin{remark}}
	\newcommand{\ermk}{\end{remark}}
\newcommand{\bhist}{\begin{histrmks}}
	\newcommand{\ehist}{\end{histrmks}}
\newcommand{\brmks}{\begin{remarks}}
	\newcommand{\ermks}{\end{remarks}}
\newcommand{\bverb}{\begin{verbatim}}
	\newcommand{\everb}{\end{verbatim}}
\newcommand{\bntn}{\begin{notation}}
	\newcommand{\entn}{\end{notation}}
\newcommand{\bfct}{\begin{fact}}
	\newcommand{\efct}{\end{fact}}
\newcommand{\bfcts}{\begin{facts}}
	\newcommand{\befcts}{\end{facts}}
\newcommand{\benum}{\begin{enumerate}}
	\newcommand{\eenum}{\end{enumerate}}
\newcommand{\bitem}{\begin{itemize}}
	\newcommand{\eitem}{\end{itemize}}

\renewcommand	{\o}		{\circ}

\renewcommand	{\epsilon}	{\varepsilon}
\renewcommand	{\a}		{\alpha}
\renewcommand	{\b}		{\beta}
\renewcommand	{\d}		{\delta}

\renewcommand	{\l}		{\lambda}

\newcommand		{\vk}		{\varkappa}

\renewcommand	{\:}		{\colon}

\newcommand		{\quotientmed}[2]	{{\raisebox{.2em}{$#1$}}\  \!\!\big/\!\!\ 
	{\raisebox{-.2em}{$#2$}}}
\newcommand		{\qquotientmed}[2]	{{\raisebox{.2em}{$#1$}} \ \big/\!\!\!\!\!\big/ \ 
	{\raisebox{-.2em}{$#2$}}}

\newcommand		{\fs}		{{\mathfrak s}}

\newcommand		{\td}		{\tilde{d}}

\newcommand		{\tA}		{\wt{A}}

\newcommand		{\tB}		{\wt{B}}

\newcommand		{\tG}		{\wt{G}}

\newcommand		{\tK}		{\wt{K}}

\newcommand		{\tS}		{\wt{S}}

\newcommand		{\wP}		{\widehat{P}}

\newcommand		{\cP}		{\widecheck{P}}
\newcommand		{\exterior}{\varbigwedge\mspace{.75mu}}
\newcommand		{\ext}	{\exterior}
\newcommand		{\varbigwedge}{{\scalerel*{\bigwedge}{I}}}
\newcommand		{\ewP}	{\ext\wP}

\newcommand		{\eP}	{\ext P}
\newcommand		{\ecP}	{\ext\cP}

\newcommand		{\eQ}	{\ext Q}
\newcommand		{\ePQ}	{\ext Q \ox \ext P}

\newcommand		{\qf}			{q^{\from}}
\DeclareMathOperator{\Kdim}	{Krull\ dim }

\newcommand		{\SSS}	{Serre spectral sequence\xspace}

\newcommand		{\eqf}	{equivariantly formal\xspace}

\newcommand		{\eqfity}	{equivariant formality\xspace}
\newcommand		{\isotf}	{isotropy-formal\xspace}
\newcommand		{\isotfity}	{isotropy-formality\xspace}

\DeclareMathOperator{\linspan}{span}

\newcommand	{\fa}{\f a}
\newcommand	{\fb}{\f b}

\newcommand{\HKt}{H_K^{\ox 2}}

\newcommand{\wHG}{\Hp_G}
\newcommand{\wHK}{\Hp_K}

\newcommand{\Hp}{H^{**}}
\newcommand{\HGp}{\Hp_G}
\newcommand{\HKp}{\Hp_K}

\newcommand{\compl}{\!\wh{\ ^{\phantom{x}}}\mn}
\newcommand{\K}{K^*}
\newcommand{\KG}{\K_G}

\newcommand{\HKGK}{\H_K G_K}

\newcommand{\wPz}{\wP_0}

\newcommand{\Gad}{G^{\Ad}}

\newcommand{\Kad}{K^{\Ad}}
\newcommand{\KGad}{\K_{\Gad}}

\newcommand{\KKad}{\K_{\Kad}}

\newcommand{\imi}{R}
\newcommand{\imd}{\im \d}

\newcommand{\keri}{\ker i^*}

\newcommand		{\Kthy}			{$K$-theory\xspace}
\newcommand		{\Ktic}			{$K$-theoretic\xspace}
\newcommand		{\Ktcly}			{$K$-theoretically\xspace}
\newcommand{\RKEF}{rationally \Ktcly \eqf}
\newcommand{\RKEFity}{rational \Ktic \eqfity}
\newcommand{\CEF}{cohomologically equivariantly formal\xspace}
\newcommand{\WEF}{weakly equivariantly formal\xspace}
\newcommand{\QWEF}{$\Q$--weakly equivariantly formal\xspace}
\newcommand{\cpair}{compact pair\xspace}

\newcommand{\ccpair}{compact, connected pair\xspace}

\newcommand{\CIR}{complete intersection ring\xspace}
\newcommand{\TFAE}{The following are equivalent\xspace}

\newcommand{\Dfrac}[2]{%
	\ooalign{%
		$\genfrac{}{}{1.2pt}0{#1}{#2}$\cr%
		$\color{white}\genfrac{}{}{.4pt}0{\phantom{#1}}{\phantom{#2}}$}%
}

\renewcommand	{\ring}		{k}

\newcommand		{\CGA}		{\textsc{cga}\xspace}
\newcommand		{\CGAs}		{\textsc{cga}s\xspace}
\newcommand		{\DGA}		{\textsc{dga}\xspace}

\newcommand		{\CDGA}		{\textsc{cdga}\xspace}
\newcommand		{\CDGAs}	{\textsc{cdga}s\xspace}

\newcommand		{\act}		{\,\raisebox{.25ex}{$\curvearrowright$}\,}

\newcommand		{\ang}[1]			{\langle #1 \rangle}

\newcommand		{\eqn}[1]			{\begin{align*} #1 \end{align*}}
\newcommand		{\quation}[1]		{\begin{equation} #1 \end{equation}}

\newcommand		{\mn}				{\mspace{-2mu}}

\newcommand		{\nd}				{\noindent}
\newcommand		{\ol}				{\overline}
\newcommand		{\os}			{\overset}
\newcommand		{\us}			{\underset}
\newcommand		{\ul}			{\underline}
\newcommand		{\wh}			{\widehat}

\newcommand		{\wt}			{\widetilde}

\newcommand		{\mr}			{\mathrm}
\newcommand		{\bb}			{\mathbb}

\newcommand		{\f}			{\mathfrak}

\newcommand		{\e}			{\epsilon}
\newcommand		{\z}			{\zeta}

\newcommand		{\s}			{\sigma}

\newcommand		{\vp}			{\varphi}

\newcommand		{\G}			{\Gamma}

\newcommand		{\pt}		{\mr{pt}}

\newcommand		{\F}		{\bb F}
\newcommand		{\N}		{\bb N}
\newcommand		{\Z}		{\bb Z}
\newcommand		{\Q}		{\bb Q}
\newcommand		{\R}		{\bb R}
\newcommand		{\C}		{\bb C}

\DeclareMathOperator{\id}		{id}

\renewcommand 		{\H}	{H^*}
\newcommand 		{\HG}	{{\H_G}}

\newcommand 		{\HK}	{{\H_K}}

\newcommand 		{\HS}	{{\H_S}}

\newcommand		{\HSN}	{H_S^N}

\newcommand		{\KGK}	{{}_K G_K}
\newcommand		{\HKN}	{H_K^N}

\let\inter\cap%

\newcommand		{\quot}		{\,/ \mn\mn /\,}

\renewcommand	{\-}		{^{-1}}
\renewcommand	{\o}		{\circ}
\renewcommand	{\.}		{\cdot}
\newcommand		{\x}		{\times}

\newcommand		{\xt}[3]	{{    {#2}\us{#1}\ox{#3}   }}

\newcommand		{\semidirect}	{\rtimes}
\newcommand		{\tensor}		{\otimes}

\newcommand		{\ox}			{\otimes}

\newcommand		{\+}			{\oplus}

\newcommand		{\Direct}		{\bigoplus}

\newcommand		{\limit}		{\varprojlim}
\newcommand		{\colim}		{\varinjlim}

\DeclareMathOperator{\diag}		{diag}

\DeclareMathOperator{\rk}		{rk }
\DeclareMathOperator{\im}		{im }

\DeclareMathOperator{\Tor}		{Tor}

\DeclareMathOperator{\ch}		{ch }

\DeclareMathOperator{\Ad}		{Ad }

\newcommand		{\U}			{\mr{U}}
\newcommand		{\SU}			{\mr{SU}}

\newcommand		{\Sp}			{\mr{Sp}}
\newcommand		{\Spin}		{\mr{Spin}}

\newcommand		{\longto} 		{\longrightarrow}
\newcommand		{\lt}			{\longto}
 
\newcommand		{\lmt}			{\longmapsto}

\newcommand		{\from}		{\leftarrow}

\newcommand		{\inc}		{\hookrightarrow}
\newcommand		{\xinc}		{\xhookrightarrow}
\newcommand		{\longinc}		{\xinc[]{\ \ \ \ }}

\newcommand		{\longmono}	{\xmono[]{\ \ \ \ }}

\newcommand		{\longepi}	{\xepi[]{\ \ \ \ }}

\newcommand		{\simto}		{\xrightarrow{\sim}}

\newcommand		{\longsimto}	{\os\sim\longto}
\newcommand		{\isoto}		{\longsimto}

\newcommand		{\ceq }			{\coloneqq}
\newcommand		{\eqc}			{\eqqcolon}
\newcommand		{\ideal}			{\unlhd}
\newcommand		{\idealneq}		{\lhd}

\newcommand		{\iso}				{\cong}
\newcommand		{\homeo}			{\approx}

\numberwithin{equation}{section}

\newcommand\blfootnote[1]{%
	\begingroup
	\renewcommand\thefootnote{}\footnote{#1}%
	\addtocounter{footnote}{-1}%
	\endgroup
}

\title{Equivariant formality of isotropy actions}
\author{Jeffrey D. Carlson\thanks{\ %
		J.D.C. was partially supported by a postdoctoral fellowship 
		from the Instituto Nacional de Matem\'{a}tica Pura e Aplicada
		(IMPA)
		funded by the 
		Coordena{\c{c}}{\~a}o de Aperfei{\c{c}}oamento 
		de Pessoal de N\'{i}vel Superior 
		(CAPES),
		and partially supported by the 
		National Center for Theoretical Sciences (Taipei) during the conclusion of writing.
	}\ \and Chi-Kwong Fok\thanks{\ C.-K. F. was supported by a postdoctoral fellowship from the National Center for Theoretical Sciences.\smallskip}}

\begin{document}

\maketitle


\begin{abstract}
	Let $G$ be a compact connected Lie group and $K$ a connected Lie subgroup.
	In this paper, we collect an assortment of results on equivariant 
	formality of the isotropy action of $K$ on $G/K$.
	If the isotropy action of $K$ on $G/K$ is equivariantly formal,
	then $G/K$ is formal in the sense of rational homotopy theory. 
	This enables us to strengthen a theorem of Shiga--Takahashi 
	to a characterization of equivariant formality in this case. 
	Using a \Ktic analogue of equivariant formality
	introduced and shown by the second-named author 
	to be equivalent to equivariant formality in the usual sense, 
	we provide a representation-theoretic characterization 
	for equivariant formality of the isotropy action 
	and give a new, uniform proof of equivariant formality 
	for some classes of homogeneous spaces for which it was
	previously known.
\end{abstract}


\section{Introduction}

Equivariant formality is an important property of 
certain topological group actions, first 
named by Goresky, Kottwitz, and MacPherson~\cite[{\SS}1.2]{GKM1998}
but already identified 
as important as far back as the Borel \emph{Seminar on transformation groups} 
\cite[Ch. XII]{borel1960seminar},
which allows, inter alia,
the application of powerful integral localization formulas~%
\cite{BV1982,AB1984,jeffreykirwan1995}.
Broad classes of actions of especial interest are well known to be \eqf, 
e.g., Hamiltonian torus actions on compact symplectic manifolds 
and linear algebraic torus actions on smooth complex projective varieties~%
\cite[\SS1.2, Thm.~14.1]{GKM1998}. 
It would be desirable to have more explicit characterizations of equivariant formality.

A natural place to start is with homogeneous spaces $G/K$,
the orbits of Lie group actions.
If $G$ is compact, the left translation action of $G$ on the right quotient $G/K$ 
is known to be \eqf if and only if $K$ contains a maximal torus of $G$; in fact, 
for any subgroup $H$ of $G$ of rank higher than that of $K$,
it is impossible that the restricted action $H$ on $G/K$ be equivariantly formal.
The next natural task, then,
and the subject of the present paper, 
is to characterize when the restricted \emph{isotropy action} 
of $K$ on $G/K$ is \eqf.
We call a pair of compact, connected Lie groups $(G, K)$ an \emph{\isotf} pair 
when this occurs.
Isotropy-formality has been studied by various authors, 
whose results we summarize here and in \Cref{sec:prev}. 
All known results involve the notion of \emph{formality} 
in the sense of rational homotopy theory. 
For instance, Shiga and Takahashi provided 
the following sufficient conditions for \isotfity.

\begin{theorem}[{\cite[Thm.~A, Prop.~4.1]{shiga1996equivariant}%
\cite[Thm.~2.2]{shigatakahashi1995}}]\label{thm:ST}
If $G/K$ is formal and $\H(BG; \Q) \lt \H(BG; \Q)^{N_G(K)}$ surjective, 
then $(G,K)$ is \isotf.
If $G/S$ is formal with $S$ a torus containing regular elements of $G$,
then $(G,S)$ is \isotf if and only if $\H(BG;\Q) \lt \H(BS; \Q)^{N_G(S)}$ 
is surjective.
\end{theorem}

Recently, Goertsches and Noshari showed that an important class of formal homogeneous spaces is \isotf.

\begin{theorem}[{\cite{goertschesnoshari2016}%
		\cite{goertsches2012isotropy}}]\label{thm:GNsurj} 
If $(G,K)$ is a (generalized) symmetric pair (see Example \ref{eg:eqf}.\ref{eg:symm}), then it is \isotf.
\end{theorem}

The first-named author of the present paper also found the following characterization 
when $K \iso S^1$.
 
\begin{theorem}[{\cite[Algorithm 1.4]{carlson2014equivariant} %
		and Example \ref{eg:eqf}.\ref{eg:TNHZ},\ref{eg:circle}}]\label{thm:reflectedcirc} 
Let $G$ be a compact, connected Lie group and $S$ a circle subgroup. 
The pair $(G, S)$ is \isotf if and only if
	\begin{enumerate}
		\item the restriction map $H^1(G; \Q) \lt H^1(S; \Q)$ is surjective, or 
		\item the restriction map $H^1(G; \Q) \lt H^1(S; \Q)$ is not surjective, and there exists $g\in G$ such that $gzg^{-1}=z^{-1}$ for all $z\in S$. 
	\end{enumerate}
These cases are mutually exclusive.
\end{theorem}

We observe that in \Cref{thm:ST}, formality is assumed, 
and the homogeneous spaces considered in 
Theorems \ref{thm:GNsurj} and \ref{thm:reflectedcirc} are formal. 
It turns out that this is no coincidence.
The first main result of this paper is as follows.

\begin{thmx}\label{thm:isotformal}
	If a pair $(G,K)$ of compact, connected Lie groups is \isotf, 
	then $G/K$ is formal.
\end{thmx}

This enables us to obtain the following 
characterization of \isotfity, 
which is a strengthened version of \Cref{thm:ST}.

\begin{restatable}{theorem}{STpp}\label{thm:ST++}
Let $(G,K)$ be a pair of compact, connected Lie groups
and $S$ a maximal torus of $K$. 
Then $(G,K)$ is \isotf 
if and only if 
it is formal and $\H(BG; \Q) \lt \H(BS; \Q)^{N_G(S)}$ is surjective.
\end{restatable}

In the latter part of this paper we appeal to another tool, namely \Kthy, to investigate \isotfity. 
Inspired by the notion of \emph{weak equivariant formality} introduced 
by Harada and Landweber~\cite{haradalandweber2007}, 
the second-named author of the present paper defined 
the related notion of 
\emph{rational K-theoretic equivariant formality} 
(\emph{RKEF} for short)~\cite{fok2017formality}, 
which amounts to surjectivity of the forgetful map 
from equivariant \Kthy with $\Q$ coefficients to ordinary \Kthy. 
He also proved the following equivalence.

\begin{restatable}[Fok~{\cite{fok2017formality}}]{theorem}{equivthm}%
\label{thm:equiv}
An action of a compact Lie group $G$ on a finite CW complex $X$
is equivariantly formal if and only if it is rational $K$-theoretic equivariantly formal,
if and only if it is \QWEF in the sense of Definition \ref{def:weakKEF}.
\end{restatable}

This equivalent formulation translates the problem of determining \isotfity 
to the context of \Kthy.
One advantage of this approach is that it becomes more straightforward to check
if the forgetful map is surjective, since this amounts to determining 
if any vector bundle on a given homogeneous space can be equipped with an equivariant structure. 
Our second main result is a representation-theoretic characterization of \isotfity.
Here, as usual, $R(\G)$ is the complex representation ring of 
a group $\G$, a ring $R$ is said to be \emph{regular at} a prime ideal $I$ if the localization $R_I$ 
is a regular local ring,
and covering homomorphism of topological groups is said to be 
\emph{central covering} is its kernel lies in the center of its domain.

\begin{thmx}\label{thm:regularity}
	Let $(G, K)$ be a compact, connected pair, 
	$\tG \lt G$ any finite, central covering 
	such that $\pi_1 \tG$ is torsion-free,
 	and $\tK$ the identity component of the preimage of $K$ under the covering map. 
 	Then $(G, K)$ is  \isotf if and only if the image $R \ceq \wt{\imath}^*R(\tG) \ox \Q$
 	of the restriction map $\wt{\imath}^*\: R(\tG) \ox \Q \lt R(\tK) \ox \Q$ 
 	is regular at the restriction $I = \wt{\imath}^*I(\tG) \ox \Q$
 	of the augmentation ideal of $R(\tG) \ox \Q$.
\end{thmx}
\Cref{thm:regularity} allows a new, uniform proof of \isotfity 
(and with \Cref{thm:isotformal}, of formality) 
for many known examples in the literature, 
including the generalized symmetric spaces of \Cref{thm:GNsurj}, 
for which the original proof by contrast requires a case-by-case analysis
running through the classification theorem for such spaces.

There is a common philosophical thread
underrunning these apparently disparate observations:
\isotfity arises in situations where \emph{symmetry} is maximized, 
as indicated by \Cref{thm:GNsurj}.
In general, the map $\H(BG;\Q) \lt \H(BS;\Q)^{N_G(S)}$ is not surjective, 
but if it is, morally speaking it is because the codomain $H(BS;\Q)^{N_G(S)}$
is minimized, meaning the size of the image of $N_G(S)$ 
in the automorphism group of $S$ is maximized.
Similarly, the generators of the kernel of $i^*\: R(G) \lt R(K)$ 
can be thought of as equations which cut out the subgroup $K$. 
The embedding of $K$ is symmetric if those equations are relatively simple, and
regularity is the condition which describes this simplicity.

The organization of the paper is as follows. 
In \Cref{sec:formality} we provide relevant definitions and notation and 
review known examples and reduction results of \isotfity, 
then recall the algebraic definitions necessary to prove 
\Cref{thm:isotformal}.
In \Cref{sec:formalrefl} we prove \Cref{thm:isotformal} 
and some natural corollaries involving invariant theory. 
In \Cref{sec:excounterex} we illustrate by way of examples 
that none of the sufficient conditions for \isotfity discussed in previous sections implies the others and thus they are not extraneous. 
In \Cref{sec:KHcompare}, we recall the second-named author's definition of \RKEFity \cite{fok2017formality}
and provide new proofs of \Cref{thm:equiv} 
and some related results from the same paper. \Cref{sec:thm:KGKstructure} concerns the rational $K$-theory ring structure of homogeneous spaces, which is parallel to the corresponding cohomological result given by \Cref{thm:purestructure}. 
Finally in \Cref{sec:regularity} we exploit \RKEFity to prove \Cref{thm:regularity} 
and give an alternative, uniform proof of \Cref{thm:GNsurj}.
Two computational examples are also given to demonstrate the utility of \Cref{thm:regularity}.  

\medskip 
 
 \nd\emph{Acknowledgments.} 
 J.D.C. gratefully acknowledges
 enlightening conversations with 
 Omar Antol\'{i}n Camarena,
 Matthias Franz,
 Oliver Goertsches, 
 Steve Halperin, 
 Friedrich Knop,
 Larry Smith, and
 Loring W.~Tu.
 C.-K.~F. would like to thank Loring W.~Tu for bringing 
 this problem to his attention, Reyer Sjamaar for his interest and Nan-Kuo Ho for her encouragement during the writing of this paper.
The authors would additionally like to thank Manuel Amann for the correction explained
in \Cref{eg:SU(7)}.

\section{Definitions and background}\label{sec:formality}
In this section we set up some notation and provide some background lemmata on equivariant cohomology and graded algebras, as well as models for homogeneous spaces and homotopy biquotients.

\subsection{Equivariant formality and previous work}\label{sec:prev}
In this subsection we summarize what is known about \isotfity,
taking the opportunity to establish some notation and terminology along the way.

\begin{definition}\label{def:eqf}
In all that follows, absent explicit indication to the contrary, 
cohomology of spaces is singular with \textbf{rational} coefficients.
Given the continuous action of a group $G$ on a space $X$,
we say the action is
\defd{equivariantly formal}
when the fiber inclusion 
in the {Borel fibration}
 $X \to (X \x EG)/G \eqc \defm{X_G} \to BG$
induces a surjection 
$\defm{\HG(X)} \longepi \H X$ from Borel equivariant cohomology
to singular cohomology,
or equivalently~\cite[Prop.~II.4.3, p.~89]{smith1967emss} if the projection 
makes $\HG X$ a free module over 
the coefficient ring $\defm{\HG} \ceq \HG(\pt) = \H(BG)$.

When $G$ is a compact Lie group and $K$ a closed subgroup, 
we say for concision that $(G,K)$ is a \defd{\cpair}.
If in addition both groups are connected, we say $(G,K)$ is a \defd{\ccpair}.
The left action of $K$ on the right quotient $G/K$, 
given by $\smash{k\. gK = (kg)K}$, 
is called the \defd{isotropy action}.
Later, we will have occation to use the observation
$EK \x_K G/K$ is homotopy equivalent to
the \defd{homotopy biquotient} $\defm\KGK$, 
the homotopy quotient of $G$ by the two-sided
$K^2$-action given by 
$(k',k)\.g = k' g k\-$.
If the isotropy action of $K$ on $G/K$ 
is equivariantly formal,
we say the pair $(G,K)$ is \defd{\isotf}.
We write $\defm{N} = \pi_0 N_G(K) = N_G(K)/Z_G(K)$ 
for the group of automorphisms of $K$ induced by conjugation,
$\defm{W_K}$ for the Weyl group of $K$,
and $\smash{\defm{\HKN} \ceq (\HK)^{N_G(K)}}$
for the invariant subring.
\end{definition}


Existing work on \isotfity
consists of an equivalent condition and several sufficient conditions.
\Cref{thm:ST}
is the most general of the sufficient conditions; 
we will improve this to an equivalence in \Cref{sec:formalrefl}.
Other results have the flavor of reductions.
For instance we can always replace $K$ by a torus
and $G$ by the product of a simply-connected group and a torus.

\bthm[Carlson {\cite[Thms.~1.1,2]{carlson2014equivariant}}]\label{thm:eqftorusreplacement}
Let $(G,K)$ be a \ccpair, 
$S$ a maximal torus of $K$,
and $\wt S$ the identity component of the preimage of $S$
under a finite covering $\wt G$ of $G$.
\bitem
\item
The pair $(G,K)$ is \isotf if and only if $(G,S)$ is.
%
\item
The pair $(G,S)$ is \isotf if and only if $(\tG, \tS)$ is. 
\eitem
\ethm

%

The following are examples of isotropy-formal pairs from the literature.

\begin{example}\label{eg:eqf}
A \cpair $(G,K)$ is known to be \isotf if any of the following conditions holds.
	\begin{enumerate}
		\item\label{eq:equal} The ranks of $K$ and $G$ are equal~%
			(Goresky et al.~%
			\cite[Thm.~14.1(1)]{GKM1998}%
			\cite[Prop.~10.3.1]{carlsonmonograph}%
			).
			
			\smallskip
			
			
		\item\label{eg:TNHZ} 
			The restriction map $\H G \lt \H K$ is surjective~%
			(Shiga \cite[Cor.~4.2]{shiga1996equivariant}%
			\cite[Prop.~10.3.2]{carlsonmonograph}).
			
			\smallskip
			
			
		\item\label{eg:symm} 
			There is a Lie group automorphism $\s$ of $G$ such that 
			$K = (G^{\s})_0$ is the identity component of the fixed point subgroup~%
			(Goertsches--Noshari~%
			\cite{goertsches2012isotropy,goertschesnoshari2016}).
			In this case, 
			we call $(G,K)$ a \defd{generalized symmetric pair}. 
			If $\s$ has finite order,
			$G/K$ is traditionally called a \emph{generalized symmetric space}, 
			and when $\s$ is an involution, 
			a \emph{symmetric space}.
			
%
		
		\item\label{eg:circle}
			If the group $K \iso S^1$ is a circle and 
			$H^1 G \lt H^1 S^1$ is \emph{not} surjective 
			(i.e., if Example \ref{eg:eqf}.\ref{eg:TNHZ} does not apply),
			then $(G,S^1)$ is \isotf if and only if $|N| = 2$, 
			meaning conjugation by some $g \in G$ induces 
			$z \mapsto z\-$ on $S^1$ 
			(Carlson \cite[Thm.~7.2]{carlson2014equivariant}). 
			
		\end{enumerate}
\end{example}

\subsection{Algebraic notions}\label{sec:pure}

To prove \Cref{thm:isotformal},
we will need to use an algebraic model of the map $G/K \lt {}_K G/K$,
so we briefly state what such models are and
basic structural results we will call on.

\bdefn\label{def:topformal}
We reserve the letters $\defm{P}$ and $\defm{Q}$ for positively-graded rational vector spaces concentrated respectively in odd and even degrees,
writing $\defm{\eQ} \ox \defm{\eP}$ 
for the free $\Q$-\CGA on this space;
$\eQ$ is a symmetric algebra and $\eP$ an exterior algebra. 
A \defd{pure Sullivan algebra}
is a $\Q$-\CDGA of the form $(\ePQ,d)$
where $d$ is a derivation such that 
\[
dP \leq \eQ \qquad \textrm{and}\qquad d(\eQ) = 0.%
\footnote{\
	The nilpotence condition~\cite[Def., p.~138]{FHT} 
	required of a Sullivan algebra
	follows automatically from 
	the restrictions on the differential.
}
\]
A \DGA $(A,d)$ is said to be \defd{formal} if it can be
joined by a zig-zag of \DGA quasi-isomorphisms
to its own cohomology $\big(\H(A,d),0\big)$, 
viewed as a \DGA with differential zero.
A $\Q$-\CDGA $(\ext V,d)$ is a \defd{model} of a topological space $X$ if it 
can be
joined by a zig-zag of \DGA quasi-isomorphisms
to the rational cochain algebra $C^*(X;\Q)$. 
A map $(\ext V_Y,d) \lt (\ext V_X,d)$ 
of models for $X$ and $Y$ is a \defd{model} of $f\: X \lt Y$ if
the induced map in cohomology is $\H f$.
The space $X$ is said to be \defd{formal} 
if $C^*(X;\Q)$ 
is formal.%
\footnote{\ 
	The standard definition asks it be weakly equivalent
	\emph{through \CDGAs} to the minimal model of $C^*(X;\Q)$,
	but recent work of Saleh~\cite{saleh2016formality}
	shows that for $\F$ of characteristic zero, 
	these notions are equivalent.
}
\edefn

Pure Sullivan algebras have a structure theory which we will need.
First note~\cite[p.~76]{koszul1950transgression} {\cite[p.~435]{FHT}} 
the exterior degree 
on $\eP$ 
induces an \defd{exterior degree} or \emph{lower grading} 
$\defm{(\eQ \ox \eP)_p} \ceq \eQ \ox \ext\mn^p P$
which persists in cohomology.
There are canonical \CDGA maps
\quation{\label{eq:CDGAsequence}
	(\eQ,0) \os \chi\lt (\ePQ,d) \os j\lt (\eP,0),
}
respectively including $\eQ \ox \Q$ in exterior degree zero and modding out $(\ext^{\geq 1} Q)$.
The image of $j^*$ is~{\cite[Thm.~13.2]{koszul1950homologie}\cite[Thm.~10.4]{andre1962tohoku}} 
the exterior algebra $\ewP$
on the \defd{Samelson space} $\defm \wP \ceq  P \cap \im j^*$.
We write $\defm{\cP} \leq P$ for a graded linear complement to $\wP$.

\bprop[{\cite[Prop. II.4.IV, p.~71]{GHVIII}\cite[pp.~141, 210, 213]{onishchik}}]%
\label{thm:explicitSamelson}
The Samelson space $\wP$ is explicitly given as
$
\{z \in P : dz \in dP \. \ext^{\geq 1}Q \},
$
verbally, those generators $z$ such that $dz$ is
redundant as a generator of $(dP)$.
The complement $d\cP$ irredundantly generates the ideal $(dP)$.
\eprop



The cochain maps $\chi$ and $j$ of (\ref{eq:CDGAsequence}) in fact determine a
factorization of cochain complexes, 
yielding the main cohomological structure theorem.

%

\bthm[{%
	\cite{andre1962tohoku}%
	\cite[pp.~73, 83, 152]{GHVIII}%
	\cite[p.~141, 211]{onishchik}%
	\cite[Thm.~7.4.7,8]{carlsonmonograph}%
}]%
\label{thm:purestructure}
Let $(\ePQ,d)$ be a pure Sullivan algebra. Then
one has an algebra decomposition 
\[
\H(\eQ \ox \eP) \iso 
\H(\eQ \ox \ecP) \ox \ewP \iso
\Big(\quotientmed{\eQ}{(dP)} \,\+\, \fa\Big) \ox\mspace{1.75mu}\ewP,
\]
where $\fa = \Direct_{p \geq 1} \H(\eQ \ox \ecP)_p$ 
is the ideal of elements of positive exterior degree.
\TFAE:
\benum
\item
The ideal $\fa$ is $0$.
\item
The algebra $(\ePQ,d)$ is formal.
\item
The map
$
(\ext Q \ox \ext P,d) 
\,\lt \,
\Big(\quotientmed{\ext Q}{(dP)} \,\ox\, \ewP,0\mspace{1mu}\Big)
$
is a quasi-isomorphism.
\item  
The ideal $(dP) \ideal \eQ$ is generated by a regular sequence,
or in other words $\eQ/(dP)$ is a complete intersection ring.
\eenum
\smallskip
If $\H(\ePQ)$ is finite-dimensional, 
one also has the equivalent condition
\smallskip
\benum
\item[(v)] $\dim P - \dim Q = \dim \wP$
(the inequality $\dim P - \dim Q \geq \dim \wP$ always holds).
\eenum
\ethm

A key lemma in the proof of \Cref{thm:isotformal}
will involve identifying when the cohomology of such a model
is free over a certain subring.
For this we will need a bit more notation.

\bdefn\label{def:indecomposables}
Given an {augmented algebra} $A \to \ring$ 
over a unital commutative ring $\defm{\ring}$,
we write 
$\defm{\tA}$ 
for the {augmentation ideal}. 
When $A = \Direct_{n \geq 0} A_n$ is a commutative $\N$-graded
algebra with $A_0 = \ring$ 
(henceforth a \defd{connected} \defm{$\ring$-\CGA}),
we only ever consider that augmentation with kernel 
$\Direct_{n \geq 1} A_n$.


Given an $A$-module $M$, 
we 
define the module of \defd{indecomposables} as the quotient
\[
\defm{Q_A}M \,\ceq \, \quotientmed{M\; }{\; \tA \.  M} \,\iso\, 
\xt A {\ring\mspace{1mu}} M,
\]
yielding a right exact functor from $A$-modules to $\ring$-modules
which if $A$ is a connected \CGA also reflects epimorphisms
of nonnegatively graded $A$-modules~%
\cite[Prop.~3.8]{milnormoore}\cite[Cor.~A.1.2]{neuselsmith}.
In case $A \lt B$ is a map of augmented $\ring$-algebras,
we also write $\defm{B \quot A} \ceq Q_A B = \ring \ox_A B$.
The quotient maps $\defm{q_{A,M}}\: M \longepi Q_A M$ 
assemble into a natural transformation $q_A\:\id \lt Q_A$.
For $M = \tA$, the indecomposables
\[
\defm{Q}A \,\ceq \, \quotientmed{\tA\; }{\; \tA \.  \tA}
\]
yield a right exact functor from augmented $\ring$-algebras to $\ring$-modules~%
\cite[Prop.~3.11]{milnormoore}.
We write $\defm q\: A \epi \wt A \epi QA$ as well
for the projection to the indecomposables.
\edefn

In case $A$ is a connected $\ring$-\CGA and $M$ is $\N$-graded, 
liftings of indecomposables are generators in the sense
that a $\ring$-linear section $\s$ of $q_{A,M}$ 
induces an $A$-module surjection 
\quation{\label{eq:generators}
	\begin{aligned}
		\defm{\psi(\s)}\:\xt \ring A {Q_A M} &\longepi 	M,\\
		a \ox \ol m\phantom{Q_A\, } &\lmt	\	a \s(\ol m),
	\end{aligned}
}
since $Q_A \psi(\s)$ is an isomorphism;
in words, $\im \s$ is a $\ring$-module of $A$-module generators.

\section{Formality and reflections}\label{sec:formalrefl}

In this section we prove \Cref{thm:isotformal} and explore its consequences.

\subsection{Construction of the model}

There are known models for $G/K$ and ${}_K G/K$
which can be constructed using standard results 
about the rational homotopy theory of pullback fibrations%
~{\cite[Props.~15.5,8]{FHT}}.
The expected algebraic map between these two 
models the fiber inclusion $G/K \longinc {}_K G/K$.

Recall Hopf's theorem~\cite[Satz I, p. 23]{hopf1941hopf} that $\H G$
is the exterior algebra on the subspace 
$\defm{P\H G}$ of primitive elements
of the coproduct on $\H G$ induced by the multiplication 
of $G$, 
and Borel's theorem~\cite[Thm.~19.1]{borelthesis} 
that the transgression $\defm\tau\: P\H G \isoto Q \HG$
in the \SSS of the universal bundle is an isomorphism
and $\HG$ can be seen as a polynomial 
algebra on a lift of $Q\HG$. 
The \defd{Cartan algebra}%
~{\cite[Thm.~5, p.~216]{cartan1950transgression}%
	\cite[Thm.~11.5.II, p.~462]{GHVIII}%
	\cite[Thm.~7.1.12]{carlsonmonograph}}\footnote{\ %
	There is also a distinct
	differential form--based \CDGA model of equivariant cohomology 
	called the \emph{Cartan model}.
	Cartan proved the result at hand using such a model, 
	in the prototypical application of equivariant cohomology.
	Borel's proof in his thesis~\cite[{\SS\SS}24--25]{borelthesis} 
	is a predecessor of the 
	rational homotopy--theoretic argument.
}
is the model 
of $G_K$ or equivalently\footnote{\ 
	It is important and not completely trivial 
	that under this substitution the maps from $G$ and to $BK$ 
	remain the expected ones up to homotopy~%
	\cite[Sec.~7.1.1]{carlsonmonograph} 
	but here we take this subtlety as dealt with.
}
of $G/K$ given by
\[
(\HK \ox \H G,d)
\] 
for $\defm d$ the unique derivation vanishing on $\HK$ 
and defined on $P\H G$ as the composition 
\[
\defm d\colon P \H G
\os\tau\longleftrightarrow 
Q\HG 
\xmono{\qf} 
\HG 
\os{\rho^*}\lt 
\HK,
\]
where
$\defm{\qf}\: Q\HG \longmono \HG$ 
is some section of the projection $q$ 
to the indecomposables
and $\rho = B(K \inc G)$.

Kapovitch~\cite[Prop.~1]{kapovitch2002biquotients}\cite[\SS3.4.2]{FOT}, %
building on work of Eschenburg~\cite{eschenburg1992biquotient},
discussed a pure Sullivan model for a
\emph{biquotient}, the orbit space
of a free two-sided action on $G$ by a closed
subgroup $U$ of $G^2$;
our part is to 
observe this construction produces a model of the
\emph{homotopy biquotient} independent of freeness
of the action. In the particular case $U = K^2$,
the resulting model is
\[
(\HK \ox \HK \ox \H G,\td),
\]
where the derivation $\defm{\td}$ 
vanishes on $\HK \ox \HK$ and
sends a primitive $z \in P\H G$ to
 \[
 \td z 
 \ceq 1 \ox dz - dz \ox 1
 \]
 for $d$ the differential in the Cartan algebra
 just discussed.
 
 By constructing a square of bundle maps
 connecting the pullback diagrams 
 inducing these two models,
 one checks that $G_K \lt \KGK$
 is modeled by reduction mod 
 $H_K^{\geq 1} \ox \HK$:
\eqn{
	(\HK \ox \HK \ox \H G,\td) &\lt (\HK \ox \H G,d),\\
	H^{\geq 1}_K \ox \HK \ox \H G&\lt 0,\\
	1 \ox x \ox z &\lmt x \ox z.
}

\begin{notation}
	For brevity, we will write $\defm{\HKt}\ceq \HK \ox \HK$.
\end{notation}

For future reference, note that the maps
inducing the Cartan model
fit into a fiber sequence
\quation{\label{eq:fibersequence}
	K \lt G \lt G/K \lt BK \lt BG,
}
thus yielding a cohomology sequence (cf.~(\ref{eq:CDGAsequence}))
\quation{\label{eq:cohomseq}
	\HG \os{\rho^*}\lt \HK \os{\chi^*}\lt \H(G/K) \os{j^*}\lt \H G \os{i^*}\lt \H K.
}

\subsection{The proof of \Cref{thm:isotformal}}

As isotropy-formality amounts to the request 
$\HKGK$ be free over $\HK$, 
the proof is largely about characterizing free modules. 
First, note that not only do liftings of indecomposables characterize free modules,
but any lifting will do.

\blem\label{thm:freesection}
Let $k$ be a commutative ring, 
$A$ an augmented $\ring$-algebra, and $M$ an $A$-module.
Then $M$ is a free $A$-module if and only if $Q_A M$ is a free $\ring$-module 
and for some section $\s$ of $q_{A,M}\: M \longepi Q_A M$,
the map $\psi(\s)$ of (\ref{eq:generators}) is an $A$-linear isomorphism.
If $A$ is a connected $\ring$-\CGA and $M$ is a \emph{finite} free $A$-module, 
then $\psi(\s)$ is an isomorphism for \emph{any} section $\s$.
\elem
\bpf
If $Q_A M \iso \ring^{\+\l}$ for some cardinal $\l$
and there exists an $A$-module isomorphism $M \iso A \ox_k Q_A M$,
then $M \iso A \ox_k k^{\+\l} \iso A^{\+\l}$
is free over $A$.
In the other direction, 
an $A$-module isomorphism $\phi\: M \isoto A^{\+ \l}$
induces a $\ring$-module isomorphism $Q_A\phi\: Q_A M \isoto \ring^{\+ \l}$,
so $Q_A M$ is free over $\ring$
and we may identify $q_{A,M}\:M \longepi Q_A M$
with $q_{A,A^{\+\l}}\: A^{\+\l} \longepi k^{\+\l}$.
Then the obvious section $\s\: \ring^{\+\l} \longmono A^{\+\l}$
makes $\psi(\s)$ an $A$-linear isomorphism.

Supposing the additional hypotheses,
for \emph{any} section $\varsigma$ of $q_{A,A^{\+\l}}$,
the map $\psi(\varsigma)$ is a surjective $A$-module map from
$A \ox_k \ring^{\+\l} \iso A^{\+\l}$ to $A^{\+\l}$,
and so is invertible~\cite[Prop.~1.2]{vasconcelos1969finitely}{\cite[Lem.~10.15.4]{stacks}}.
\epf

Second, base extensions along connected \CGAs reflect freeness.

\begin{lemma}\label{thm:isotformalcommalg}
	Let $k$ a commutative ring, 
	$A \leq B$ connected $\ring$-\CGAs, 
	and $M$ a finite $A$-module.
	If $\smash{\xt A B M}$ is free over $B$,
	then $M$ is free over $A$.
\end{lemma}
\bpf
Applying the right exact functor $- \ox_A M$ 
to the exact sequence $0 \to \tB \to B \to \ring \to 0$ of $A$-modules,
one finds that $Q_B(B \ox_A M) \iso \ring \ox_A M$.
From the assumption that $B \ox_A M$ is free over $B$,
it follows that $\ring \ox_A M$ is free and finite over $\ring$, 
say on the basis $(1 \ox x_j)$
for some elements $x_j \in M$,
and since this basis is finite,
it follows by \Cref{thm:freesection}
that extending the structure map $\ring \to B$ 
to a section $\s\: \ring \ox_A M \longmono B \ox_A M$
of the projection $q_{B,B\ox_A M}$ 
induces a $B$-module isomorphism\footnote{\ %
	N.B. the dependence of this map on the arbitrary choice of basis;
	the ``natural'' guess taking $b \ox (1\ox x) \mapsto b \ox x$
	for all $b$ and $x$ is ill-defined.
}
\eqn{
	\xt \ring B {(\xt A\ring M)}		&\os{\defm{\psi}}\lt \xt ABM,\\
	\sum b_j \ox (1 \ox x_j)	&\lmt	\sum b_j  \ox x_j.
}
The restriction of $\psi$ to $A \ox_\ring (k \ox_A M)$ 
factors through $A \ox_A M$
as the surjective map in (\ref{eq:generators}):
\[
\xymatrix@C=1.5em@R=3.75em{
	\xt \ring A{(\xt A kM)}\ar@{^{(}->}[d] \ar[r]^(.58){\defm{\varphi}}
	&	\xt A A M\ar[r]^(.6)\sim \ar[d]
	&	M 
	\\ 
	\xt \ring B {(\xt A kM)} \ar[r]^(.58)\sim_(.58)\psi	
	& 	\xt ABM.
}
\]
This map $\varphi$ must also be injective because the composition
$A \ox_k k \ox_A M \lt B \ox_A M$
along the lower-left is.
The composition of $\varphi$ with the standard $A$-module isomorphism
$A \ox_A M \isoto M$ presents $M$ as the free $A$-module on the basis $(x_j)$.
\epf

We now apply these results to our model.

%


\begin{proof}[Proof of \Cref{thm:isotformal}]
Recall that
{\isotfity} is the two equivalent demands 
that the fiber inclusion in the Borel fibration 
$
	G_K \to \KGK \to BK
$
induce a surjection $\HKGK \lt \H(G_K)$ in cohomology
and the projection 
make $\HKGK$ a free $\HK$-module.
It follows from surjectivity that the images in $\H G$ of the maps 
$\H\big(\HKt \ox \H G\big) \to \H(\HK \ox \H G) \to \H G$
induced by the \CDGA projections
$\HKt \ox \H G \to \HK \ox \H G \to \H G$
are equal,
so the Samelson spaces of the two algebras agree.
By \Cref{thm:explicitSamelson},
$d$ and $\td$ respectively take a homogeneous basis of 
a linear complement ${\cP} \leq P = P\H G$ 
to the Samelson space $\wP$
to irredundant sets of generators for $(dP) = (d\cP) \ideal \HK$ 
and $(\td P) = (\td \cP) \ideal \HKt$.
Write $\defm A \leq \HK$ for the subring generated by 
these chosen irredundant generators $d\cP$ of $(dP) \ideal \HK$. 
The component of $\H\big(\HKt \ox \H G\big)$ of exterior degree zero is 
\[
\quotientmed{\HK \ox \HK\;}{\,(\td P)} 
\,=\, 
\quotientmed{\HK \ox \HK\;}{\,(\td \cP)}
\,=\,
\xt{A}{\HK}{\HK}
\]
since $(\td P) = (\td \cP)$ 
and $\td z = 1 \ox dz - dz \ox 1$ for $z \in P$.
Since the left $\HK$ factor of $\HKt \ox \H G$ lies in 
exterior grade zero, the left multiplication action of $\HK$ on $\H\big(\HKt \ox \H G\big)$
preserves exterior degree,
so as
$\H\big(\HKt \ox \H G\big)$ is a free $\HK$-module,
its zero-graded component $\HK \ox_A \HK$ is as well. 
Applying Lemma \ref{thm:isotformalcommalg} in the case $k = \Q$ and $B = M = \HK$, 
we see the polynomial ring $\smash{\HK}$ is free over $A$,
and thus by a classical theorem 
of Macaulay~\cite[Cor.~6.4.4]{smithinvariantbook},
the ring $\smash{A = \Q[d\cP]}$ is itself polynomial
and a basis of $d\cP$ forms a regular sequence in $\HK$.
Thus $\HK \quot \HG$ is a complete intersection ring,
so $G/K$ is formal by Theorem \ref{thm:purestructure}.
\end{proof}

\subsection{Consequences of formality}

Isotropy-formality implies formality of the homotopy quotient as well.

\bcor\label{thm:hmtquotformal}
Let $(G,K)$ be an \isotf \ccpair.
Then 
$\KGK$ is itself formal.
Particularly, 
	there is an $\HKt$-algebra isomorphism 
	\[
	\HK(G/K) \iso \big(\xt{\HG}{\HK}{\HK}\big) \ox \ewP,
	\]
	where $\ewP$ is the Samelson ring $\im\!\big(
	\mspace{-1mu}\H(G/K) \lt \H G\big)$.
\ecor
\bpf
This ring structure was proven in the first-named 
author's thesis under an additional condition
$G/K$ be formal~\cite[Thm.~11.1.1]{carlsonmonograph}
(and later published~\cite{carlson2016grassmannian})
but by \Cref{thm:isotformal},
this is already the case if $(G,K)$ is \isotf.
By \Cref{thm:purestructure}, this implies $\KGK$
is formal.
\epf

\bex\label{eg:SU(3)2}
The converse to \Cref{thm:hmtquotformal}
does not hold.
For example, consider 
the block-diagonal inclusion of $\SU(3)^2$ in  $\SU(6)$.
In the Kapovitch model for $\smash{{}_{\SU(3)^2} \SU(6)_{\SU(3)^2}}$,
one can show the five differentials of a basis of $P\H\SU(6)$
form a regular sequence in $H_{\SU(3)^2}^{\ox 2}$,
so by \Cref{thm:purestructure}(iv),
the homotopy quotient ${}_{\SU(3)^2} \SU(6)_{\SU(3)^2}$
is formal.
However, 
the differentials 
of the same primitives in 
the Cartan algebra for $\SU(6)/\SU(3)^2$
are known not to form a regular sequence,
so again by \Cref{thm:purestructure}(iv),
$\SU(6)/\SU(3)^2$ is not formal%
~\cite[App.~A]{amann2013nonformal}%
\cite[p.~486--488]{GHVIII}
and hence by \Cref{thm:isotformal}
the pair  $\big(\SU(6),\SU(3)^2\big)$
is not \isotf.
\eex

\begin{remark}\label{eg:SU(7)}
A theorem of Manuel Amann [personal communication] implies that when $(G,K)$ is a \ccpair,
then its isotropy-formality is in fact equivalent to simultaneous formality of both $G/K$ and ${}_K G_K$.
This contradicts Example 3.6 of the published version of this paper,
which rested on an erroneous calculation by the first-named author.
This remark is meant to rectify the error while preserving the theorem numbering of the published version.
\end{remark}

\bcor\label{thm:dimimageSamelson}
	A \ccpair $(G,K)$ is \isotf if and only if 
	the image of $\defm{g}\: \HKGK \to \H(G/K) \to \H G$ 
	meets $P\H G$ in a space of dimension $\rk G - \rk K$.
\ecor
\bpf
	If $(G,K)$ is \isotf,
	then it is formal by \Cref{thm:isotformal},
	and so by \Cref{thm:purestructure}
	we have $\smash{\H(G/K) \iso (\mn\HK \quot \HG) \ox \ewP}$
	with $\dim \wP = \rk G - \rk K$.
	Then since $\HKGK \lt \H(G/K)$ is surjective,
	its image particularly contains $\wP$,
	so $\im g$ does as well.
	In the other direction,
	$\im g$ lies in $\ewP$, 
	so if it meets $P \H G$ in a space of dimension $\rk G - \rk K$,
	then $\dim \wP=\rk G-\rk K$, meaning $G/K$ is formal, 
	and moreover, $\im g$ contains $\ewP$.
As $\HK \ox_{\HG} \HK \lt \HK\quot \HG$ 
as always surjective, 
the map $\HKGK \lt \H(G/K)$ is surjective as well and hence 
	$(G, K)$ is isotropy-formal.
\epf

In particular, \isotfity of $(G,K)$ is equivalent to a 
statement about a certain privileged
subspace of the Samelson space.

\bdefn\label{def:wPz}
The transgression $\tau\: P \H G \isoto Q \HG$
in the \SSS of the universal bundle 
$EG \to BG$
induces a well-defined inverse function
$\smash{\defm{\s}\: \HG \xepi{q} Q\HG \simto P\H G \inc \H G}$,
the \defd{suspension}.
Here $q$ is as defined in Definition \ref{def:indecomposables}.
We set $\defm \wPz \ceq \s \ker (\rho^*\: \HG \to \HK) \leq P \H G$.
\edefn

It is clear from \Cref{thm:explicitSamelson} that $\wPz$ 
is contained in the Samelson space $\wP$ of the Cartan algebra,
for since $\qf$ is a section of $q$,
for all $z \in \wPz$ we have
$d z = \rho^* \qf \tau z 
\in \rho^*(\ker \rho^* + \ker q)
= \rho^* H_G^{\geq 1}\.\rho^* H_G^{\geq 1}$
Note that the Cartan algebra and Kapovitch model 
compute $\HK(G/K) \lt \H(G/K)$
regardless of the choice of section $\qf$
we employ in the definition of the differentials,
so we are to choose this section at will.

\blem\label{thm:wPzimage}
The space $\wPz$ lies in the image of $g\:\HKGK \lt \H G$.
\elem
\bpf
Choose $\qf$ such that $\qf \tau \wPz = \qf q (\ker \rho^*)$
is contained in $\ker \rho^*$ on the nose, 
rather than merely modulo the decomposables $H^{\geq 1}_G\. H^{\geq 1}_G$.
Then, since the Cartan differential $d$ is 
$\rho^* \qf \tau$ on $P\H G$, we have
$
d \s \ker \rho^* = 
\rho^* \qf q \ker \rho^* \leq 
\rho^* (\ker \rho^*) = 
0,	
$
so that for $z \in \wPz$ 
we have $\td z = 1 \ox dz - dz \ox 1 = 0$
in the Kapovitch model $(\HKt \ox \H G, \td)$.
\epf

\bprop\label{thm:HSamelsondeltaformal}
A \ccpair $(G,K)$ is \isotf 
if and only if $\dim q(\ker \rho^*) = \dim \wPz = \rk G - \rk K$.
\eprop
\bpf
By \Cref{thm:wPzimage},
if $\dim \wPz = \rk G - \rk K$,
then $(G,K)$ is \isotf by \Cref{thm:dimimageSamelson}.
On the other hand if $(G,K)$ is \isotf,
then by the enhanced Shiga--Takahashi \cref{thm:ST++} 
and \Cref{thm:formalrefl} to follow,
$\defm\vk\:\HG \longepi \HSN$ 
is surjective and $\HSN$ is a polynomial ring
on the same number of indeterminates as generate $\HS$ and $\HK$,
the vector space $\HK\quot\HSN$ being finite-dimensional.
It follows 
\[
\dim \ker Q\vk = \dim Q\HG - \dim Q\HSN = \rk G - \rk K.
\]
But $\ker Q\vk = q(\ker \vk)$:\ 
an element of the former is $q(x)$ for some $x \in \HG$
such that 
$
	\vk(x) \in (H_S^{\geq 1})^N\.(H_S^{\geq 1})^N = 
	\vk(H_G^{\geq 1}\. H_G^{\geq 1}),
$
meaning there is $y \in H_G^{\geq 1}\. H_G^{\geq 1}$
such that $x-y \in \ker \vk$
and hence $q(x) = q(x-y) \in q(\ker \vk)$.
Finally, $\rho^*$ factors as $\smash{\HG \xepi[\vphantom{f}\vk]{} 
	\HSN \mono \HK}$
since $W_K \leq N$,
so $\ker \vk = \ker \rho^*$.
\epf

Later, in \Cref{sec:regularity},
we will find a \Ktic analogue 
of Proposition \ref{thm:HSamelsondeltaformal},
namely Proposition \ref{thm:KSamelsondeltaformal},
naturally constructing a class of vector bundles 
which admit equivariant lifts.
As formality is one of the conditions of the Shiga--Takahashi  \cref{thm:ST},
our \Cref{thm:isotformal} allows one to substantially strengthen it.

\STpp*
\bpf
	The ``if'' direction is Shiga's original result.
	The other direction follows from the Shiga--Takahashi theorem 
	and \Cref{thm:isotformal}.
	The regular element hypothesis\footnote{\ 
		A regular element of a Lie group $G$ is one 
		lying in a unique maximal torus.}
	turns out not to play
	an essential role in the proof of 
	\Cref{thm:ST} and can be omitted.
\epf

\brmk
	\Cref{thm:ST++} can actually be strengthened to replace surjectivity of
	$\smash{\HG \lt H_S^{N_G(S)}}$ with that of
	$\smash{\HG \lt H_K^{N_G(K)}}$, 
	but doing so would take us far afield,
	as it seems to require an alternate proof of the 
	original Shiga--Takahashi theorem.
\ermk

We also have an invariant-theoretic formulation.
%
%
%
Recall that we write $N = \pi_0 N_G(S)$.

\bprop\label{thm:formalrefl}
Suppose a \ccpair $(G,K)$ is such that $\smash{\HG \lt H_S^N}$ is surjective.
Then $(G,K)$ is \isotf if and only if it is formal 
and $\pi_0 N_G(S)$ acts on the Lie algebra $\defm{\fs}$ of $S$ 
as a reflection group.
\eprop
\bpf
Note that 
$
	\HS\quot\HG = \HS\quot\HSN
$
by the assumption $\HG \longepi \HSN$.
By \Cref{thm:ST++},
$(G,K)$ is \isotf if and only if $G/K$ is formal,
which by \Cref{thm:purestructure}
occurs if and only if $\HS \quot \HG = \HS \quot \HSN$ 
is a complete intersection ring, 
which by the Chevalley--Shepherd--Todd theorem~{\cite[p. 192]{kane}}
occurs if and only if $N$ is a reflection group.
\epf
%

\section{Examples and counterexamples}\label{sec:excounterex}

This section is devoted to showing the irredundancy of the three 
conditions for \isotfity discussed
in the previous subsection.
These non-formal examples are mostly to be found in a section 
in Greub et al.~%
\cite[Ch. XI, \SS 5]{GHVIII}.\footnote{\ 
See the paper of Amann~\cite{amann2013nonformal} for many more,
and Onishchik \cite[{\SS}13.4]{onishchik} for another family.
Note these examples are all of deficiency $\dim P - \dim Q - \dim \wP = 1$.}
Example \ref{eg:Sp(5)} for $n=5$ is due to Borel and 
according to Paul Baum has been circulating since at least the 1960s.

\bex
That $N$ be a reflection group 
does not ensure $G/K$ be formal
nor $\HG \lt \HSN$ surjective.
\smallskip

The pair $\big(\U(5),S\big)$, 
where $S$ is the four-dimensional subtorus 
$\big\{\!\diag(z^4,w^3,\z^2,zw\z, \vartheta)\big\}$ 
of the diagonal torus,
has $N = 1$ a reflection group (generated by zero reflections).
A computation with the computer algebra system Macaulay2
shows the regular sequence condition of
\Cref{thm:purestructure}(iv) is violated,
so $\U(5)/S$ is not formal.
One also incidentally sees that $\dim_\Q\HS\quot\HG = 22$,
meaning particularly that the map $\HG \lt \HSN = \HS$ is not surjective.
The relevant code is in an auxiliary file hosted on the first-named
author's website~\cite{carlsonmacaulay}.
\eex

\bex\label{eg:SU(3)/S}
That $(G,K)$ be formal and $N$ a reflection group (even the Weyl group of $K$)
does not ensure $\HG \lt \HSN$ be surjective.

\smallskip

The pair 
$\big(\SU(7),\SU(3) \x \SU(4)\big)$ from Example \ref{eg:SU(7)}
is formal but not isotropy-formal, 
and $N = S_3 \x S_4 = W_K$ is a reflection group.
As another example,
consider a \ccpair $(G,S)$, 
where $S \iso S^1$ is a circle 
not reflected by the larger group,
such as
$\big\{\!\diag(z,z,z^{-2})\big\}$ 
or $\big\{\!\diag(z,z^2,z^{-3})\big\}$ in $\SU(3)$. 
Then $N = 1$ 
and $(G,S)$ is formal~\cite[App.~A]{carlson2014equivariant},
but by Example \ref{eg:eqf}.\ref{eg:circle}, the pair $(G,S)$ is not \isotf.

\eex

\begin{example}\label{eg:SU(3)2refl}
That $\HG \lt \HKN$ be surjective
does not ensure $N$ be a reflection group or $(G,K)$ formal. 

\smallskip

Consider again the block-diagonal inclusion of $K = \SU(3)^2$ in $G = \SU(6)$ from Example \ref{eg:SU(3)2}.
We already saw $G/K$ is not formal.
The Weyl group $W_G = S_6$
permutes the six coordinates of the diagonal maximal torus $T$
and the stabilizer $N$ of $S = T \cap K$ in $W_G$
is generated by $W_K = S_3 \x S_3$ and
$\e = (1\, 4)(2\, 5)(3\, 6)$,
which is not a product of reflections of $\fs$.
But $\HG \lt \HSN = {H_K^{\ang \e}}$
is indeed surjective;
explicitly, it is
\[
\Q[c_2,c_3,c_4,c_5,c_6] \lt
{\big(\Q[c_2,c_3] \ox \Q[c'_2,c'_3]\big)^{\ang \e}}
=\Q[c_2+c'_2,\ c_3+c'_3,\ c_2c'_2,\ c_2 c_3' + c_2' c_3,\ c_3c'_3],\]
each displayed generator on the left being sent to the 
corresponding generator on the right.
\end{example}

\begin{example}\label{eg:Sp(5)}
None of the conditions need hold at all.

\smallskip

Consider the natural embedding of $K = \SU(n)$ in $G = \Sp(n)$.
The Weyl group $W_G = \{\pm 1\}^n \semidirect S_n$ 
acts on the diagonal maximal torus $T$ of $G$
by permutating and inverting coordinates,
and the condition $\det = 1$ on $S = T \cap K$ yields
$
	N = S_n \x \ang{\e},
$
where $\e(t_1,\ldots,t_n) = (t_1\-,\ldots,t_n\-)$.
It can be shown $\Sp(n)/\SU(n)$ is formal%
~{\cite[\SS11.15, pp.~488--90]{GHVIII}} 
(in fact, \isotf) if and only if $n \leq 4$
and that $N$ is a reflection group in exactly those cases. 
The map $\smash{\H_{\Sp(n)} \lt H_{\SU(n)}^{\ang\e}}$ 
on the other hand is surjective 
if and only if $n \leq 5$.
Thus $\big(\Sp(5),\SU(5)\big)$ shows, like Example \ref{eg:SU(3)2refl},
that $\HG \lt \HSN$ can be surjective without $N$ being
a reflection group or $G/S$ formal.
The larger pairs $\big(\Sp(n),\SU(n)\big)$ provide examples in which \emph{none}
of the three conditions 
figuring in \Cref{thm:formalrefl} hold.
\end{example}

In summary,\bitem
\item None of the three conditions need hold.
\item None of the three conditions alone implies \isotfity.
\item That $N$ be a reflection group implies neither of the other conditions.
\item That $\HG \lt \HSN$ be surjective implies neither of the other conditions.
\item That $(G,K)$ be formal does not imply that $\HG \lt \HSN$ be surjective.
\eitem
%
%
%
%
%

\section{Comparing equivariant \textit{K}-theory 
and cohomology}\label{sec:KHcompare}

In this section, we switch to another tool, namely (equivariant) $K$-theory, and develop $K$-theoretic results necessary to study the equivariant formality of homogeneous spaces later on. 

\begin{definition}[{\cite[Def. 4.1]{haradalandweber2007}}]
\label{def:weakKEF} 
Let $\ring$ be a torsion-free ring.
Write
$\defm{\K(-;\ring)}$ and $\defm{\KG(-;\ring)}$ for the unique additive 
($G$-equivariant) cohomology theories extending
to all \mbox{($G$--)CW} complexes the functors defined by
$X \lmt \K X \ox \ring$ and $(G \act X) \lmt \KG(X) \ox \ring$ 
on finite complexes.
We also write $\defm{R(G;\ring)} = R(G) \ox \ring$ and $\defm{I(G;\ring)} = I(G) \ox \ring$
for the extended representation ring and augmentation ideal.
The forgetful map
$
	K^*_G(X)\lt K^*(X)
$
factors through the homomorphism
\[
	\xt{R(G;\ring)}{\ring\mn}{\mn K_G^*(X; \ring)} \lt K^*(X; \ring).
\]
We say the action is $\defm{\ring}$\defd{--weakly equivariantly formal}
if the latter map is an isomorphism.
Following Harada--Landweber, we simply say the action is \defd{\WEF}
in the case $\ring = \Z$.
\end{definition}

As we defined equivariant formality as the
surjectivity of $\HG X \lt \H X$,
for us the more natural analogue in \Kthy is the following
definition due to the second-named author.

\bdefn[Fok {\cite{fok2017formality}}]\label{def:rkef}
A $G$-action on a space $X$
is said to be \defd{rational K-theoretic equivariantly formal} (or just \emph{RKEF}) if the forgetful map 
\[
	f \:\KG(X;\Q) \lt \K(X;\Q) 
\] 
is surjective.
\edefn
 
This condition admits a natural interpretation in terms of vector bundles: 
for every vector bundle $V$ over $X$ 
or its suspension $\Sigma X$,
there are natural numbers $m,n$ 
such that $V^{\+ m} \+ \ul{\C}^n$ admits an equivariant $G$-structure.

%
%
%
%

\subsection{The equivariant Chern character}

Recall \cite[\SS2.4]{atiyahhirzebruch}
that the {Chern character}
induces
a natural $\Z/2$-graded ring homomorphism 
$\smash{
	\defm{\ch}\:  \K (X) \lt \defm{\Hp X} \ceq \prod_n H^n X
}$
which becomes an isomorphism 
$\smash{
	\K(X;\Q) \isoto  \H X
}$
on finite CW complexes $X$. 
Analogously, a $\smash G$-equivariant vector bundle $\smash{V \to X}$
induces a vector bundle $V_G \to X_G$ of homotopy quotients
and a class $\defm{\ch_G}(V) \ceq \ch(V_G) \in \Hp X_G \eqc \defm{\HGp X}$,
the \defd{equivariant Chern character} of $V$.
As homotopy quotients respect the semiring operations, 
these classes collate into a natural $\Z/2$-graded ring homomorphism
\[
	{\ch_G}\:	\KG (X) \lt \HGp X
\]
which, though typically far from surjective, 
nevertheless prescribes $\HGp X$.

\bthm\label{thm:eqcherniso}
Let $X$ be a compact $G$-space such that 
$\KG X$ is a finite $R(G)$-module.
The equivariant Chern character induces $\Z/2$-graded ring 
isomorphisms 
\[
	\KG(X;\Q)\compl 
		\,\isoto\, 
	\HGp X 
		\,\xleftarrow{\, \sim\, }\, 
	\xt{R(G)}{\HGp\mn}{\mn\KG X} 
\]

\vspace{-.75em}

\nd natural in such $X$,
where $\defm{\KG(X;\Q)\compl}$
is the completion of $\KG(X;\Q)$
with respect to the augmentation ideal $I(G;\Q) \idealneq R(G;\Q)$.
This isomorphism preserves the augmentation to $\Q$.
\ethm

\bpf
Let $\defm{E_n G \to B_n G}$ be 
the restrictions of $EG \to BG$ over CW $n$-skeleta of $BG$
and $\defm{X_{n,G}}$ the compact spaces $E_n G \x_G X$.
The equivariant
Chern character $\KG (X) \to \K (X_G) \to \Hp X_G$
induces pro-ring maps
$\smash{
\big(\mspace{-1mu}\KG (X)/I(G)^n\.\KG (X)\mspace{-1mu}\big) \to 
\big(\mspace{-1mu}\K (X_{n,G})\mspace{-1mu}\big) \os\ch\to
(\H X_{n,G})}$.
By the Atiyah--Segal completion theorem 
\cite[Cor.~2.1]{atiyahsegalcompletion},
the second map is a pro-ring isomorphism,
and since all objects are finite $R(G)$-modules
and the $X_{n,G}$ are homotopy equivalent to finite CW complexes,
tensoring with $\Q$ yields pro-ring isomorphisms
\[
\big(\KG(X;\Q)/I(G;\Q)^n\.\KG(X;\Q)\big) 
        \simto 
    \big(\K(X_{n,G};\Q)\big) 
	    \simto
	(\H X_{n,G})
	.
\]
Since the first inverse system 
clearly satisfies the Mittag--Leffler condition
and thus has trivial $\limit\!^1$, so do the other two,
and hence~\cite[Lem. 2]{milnor1962axiomatic} 
taking limits yields isomorphisms
$
    \KG (X;\Q)\compl \simto \K(X_G;\Q) \simto \Hp X_G.
$
Particularly, $R(G;\Q)\compl \iso \HGp$,
and since $R(G)$ is Noetherian, 
by the finiteness assumption it follows~\cite[Prop.~10.13]{atiyahmacdonald}  
    \[
        \HGp X 
            \,\iso\,  
        \xt{R(G;\Q)}{\HGp\mn}{\mn\KG(X;\Q)}
            \,\iso\,
        \xt{R(G)}{\HGp\mn}{\mn\KG X}.
    \]
Naturality follows from cellular approximation
and naturality of the Chern character.
\epf

\brmk\label{rmk:Chern}
%
%
The first isomorphism in \Cref{thm:eqcherniso} 
has been stated in the mathematical physics literature 
and cited as ``a \emph{completion theorem} of Atiyah and Segal''%
~\cite[Thm.~6.7]{valentino2008ramond},
though technically it is a distinct corollary
requiring some argument like the above.
\ermk

We need one more commutative algebra lemma 
before the main result of this section.

\blem\label{thm:tensorreduction}
	Let $A$ be a Noetherian ring,
	$M$ a finite $A$-module,
	$\wh A$ and $\wh M$ the respective completions of 
	$A$ and $M$ with respect to an ideal $\fa$ of $A$, and 
	$N$ an $\wh A$-module.
	Then the induced $A$-module structure on $N$ is 
	such that
	$\smash{\xt A N M \iso \xt{\wh A} N {\wh M}}$.
\elem
\bpf
The finiteness assumptions make the natural map
$\smash{\wh A \ox_A M \lt \wh M}$
a $\wh A$-module isomorphism~\cite[Prop.~10.13]{atiyahmacdonald}.
Thence
$
	\xt {\wh A} N {\wh M}		 	\iso 
	\xt {\wh A} N {\xt A{\wh A} M} 	\iso 
	\xt A N M.
$
\epf
%
%
%
%

The main result of this section is the following result 
of the second-named author,
published with a different proof 
in earlier work.

\equivthm*
\bpf
Consider the commutative diagram
\[        
    \xymatrix@R=2.75em@C=1.5em{        
      \KG (X;\Q)    \ar@{->>}[r] \ar[d]^{\ch_G} 
    & \Dfrac{\KG(X;\Q)}{R(G;\Q)} \ar[r]_(.5){\defm{\bar f}} \ar[d]
    & \K (X;\Q)\ar[d]^{\ch}_[@!-90]\sim
    \\ \HGp X    \ar@{->>}[r]                
    & \qquotientmed{\HGp X}{\HGp} \ar[r]^(.6){\defm{\bar g}}
    & \H X
    }
\]
factoring the forgetful map $f\: \KG(X;\Q) \lt \K(X;\Q)$
and restriction $g\: \HGp X \lt \H X$.
By \Cref{thm:tensorreduction},
applied to the case where $A = R(G;\Q)$,\ \,$M = \KG(X;\Q)$, 
and $N = \Q$,
and the isomorphism $R(G;\Q)\compl \iso \HGp X$ 
from \Cref{thm:eqcherniso},
the middle vertical map is an isomorphism.
Thus we can identify $\bar f$ with $\bar g$.
But the \SSS of $X \to X_G \to BG$
collapses at $E_2$ if and only if $g$ or equivalently $\bar g$ is surjective,
implying the kernel of $g$ is the ideal generated by $H^{\geq 1}_G$,
so that $\bar g$ is an isomorphism.
Thus 
\[
	\bar f \mbox{ is surjective}
		\iff 
	\bar g \mbox{ is surjective}
			\iff
	\bar g \mbox{ is an isomorphism}
		\iff
	\bar f \mbox{ is an isomorphism.}
\]
By definition,
the action is 
\CEF if and only if $g$ or equivalently $\bar g$ is surjective,
\RKEF if and only if $f$ or equivalently $\bar f$ is surjective, and 
\QWEF if and only if $\bar f$ is an isomorphism.
\epf

%
%
%
%
%
%
%
%

\section{The \textit{K}-theory 
of compact homogeneous spaces}\label{sec:thm:KGKstructure}

In this section, we assume that $G$ and $K$ are compact, 
connected Lie groups unless otherwise specified.
The structure \cref{thm:purestructure},
as applied to the Cartan algebra $(\HK \ox \H G,d)$
computing $\H(G/K)$,
carries 
over to the context of \Kthy by means of the Chern character. 

For the right action of $K$ on $G$, 
the structure map 
$\a\mspace{-1mu}\: R(K) \to K^0_K G \os\sim\to K^0(G/K)$ 
making $K^*_K (G)$ an $R(K)$-algebra
sends a $K$-representation $\rho$ 
to the class of the associated vector bundle $G\times_K V_\rho$. 
Recall that $K^*(G; \Q)$ is a Hopf algebra 
with comultiplication induced by the group multiplication,
so it makes sense to speak of its space of primitives
and we can define a $K$-theoretic Samelson space in analogy with
the definition in \Cref{sec:pure}.


\begin{definition}
	The 
	\defd{$K$-theoretic Samelson space} of $G/K$, denoted by $\defm{\wP}$,
	is the space $PK^*(G; \Q) \cap \im j^*$ 
	of primitive elements of $\K(G;\Q)$
	lying in the image of the pullback 
	$j^*\: K^*(G/K; \Q) \lt K^*(G; \Q)$.
	We also use $\wP$ to denote a preimage in $K^*(G/K; \Q)$ under the map $j^*$ 
	if there is no danger of confusion, 
	given that
	$j^*$ maps isomorphically from this preimage onto $\wP$. 
	We call $\ext \wP \leq \K (G;\Q)$ (or $\ext \wP\leq \K(G/K; \Q)$) the
	$\defm K$-\defd{theoretic Samelson ring}. 
\end{definition}

\bthm\label{thm:KGKstructure}
Let $(G,K)$ be a \ccpair. 
The maps of (\ref{eq:cohomseqmap}) below induce a ring isomorphism
		\[
			\K(G/K;\Q)
				\,\iso\, 
			\Big(\qquotientmed{R(K;\Q)}{R(G;\Q)}  \,\+ \, \fa\Big)  \mspace{1mu}\ox\, \ewP.
		\]
Here ${R(K)}\quot{R(G)}$ is the image of $\a\: R(K) \lt \K(G/K)$,
the subalgebra $\ewP$ 
is taken isomorphically onto the 
image of 
		$j^*\:\K(G/K;\Q) \lt \K(G;\Q)$,
and the summand $\mathfrak{a}$ is an ideal of the other tensor factor.
We always have $\dim \wP\leq \rk G-\rk K$.
The space $G/K$ is formal 
	if and only if $\dim \wP = \rk G - \rk K$,
	if and only if $\fa = 0$, 
	if and only if $R(K;\Q)\quot R(G;\Q)$ is a \CIR.
\ethm

The Chern character isomorphism as applied to \Cref{thm:purestructure}
in the special case of the Cartan algebra $(\HK \ox \H G,d)$
already yields a decomposition of the isomorphism type indicated,
so the new content is only that the three factors can be identified as claimed.
The proof uses a commutative diagram
\quation{\label{eq:cohomseqmap}
	\begin{aligned}
		\xymatrix@R=3.75em@C=1.25em{
			R(G;\Q)		\ar[r]^{i^*}		
										\ar@{ >->}[d]^{{\ch_G}}	&
			R(K;\Q)		\ar[r]^(.425)\a		
										\ar@{ >->}[d]^{{\ch_K}}	&
			\K(G/K;\Q)	\ar[r]^(.57){j^*}		\ar[d]_[@!-90]{\sim}^{\ch} 	&
			\K (G;\Q)	\ar[r]^{i^*}			\ar[d]_[@!-90]{\sim}^{\ch}	&
			\K (K;\Q) 							\ar[d]_[@!-90]{\sim}^{\ch}	\\
			\wHG		\ar[r]_{\rho^* = (Bi)^*}							& 
			\wHK		\ar[r]_(.45){\chi^*}								&
			\H(G/K)		\ar[r]_(.53){j^*}									&
			\H G  		\ar[r]_(.54){i^*}									&	
			\H K
	}
	\end{aligned}
}
induced by (\ref{eq:fibersequence}) and expanding (\ref{eq:cohomseq}).

\bpf
To see that the image of $\a$ is $R(K;\Q)\quot R(G;\Q)$,
we show that ${\ch} \o \a$
induces an isomorphism
$\Q \ox_{R(G;\Q)} R(K;\Q) \isoto \Q \ox_{\HGp} \HKp$.
By \Cref{thm:eqcherniso},
$\HKp$ is the completion of $R(K;\Q)$ with respect to $I(K;\Q)$,
and since the $I(K)$-adic and $I(G)$-adic topologies on $R(K)$
agree~\cite[Cor.~3.9]{segal1968representation},
$\HKp$ is also the completion with respect to $i^* I(G;\Q)$.
Thus \Cref{thm:tensorreduction}
applied in the case $A = R(G;\Q)$ and $M = R(K;\Q)$ and $N = \Q$
yields the isomorphism.

The identification of the exterior factor
follows from 
the third square in (\ref{eq:cohomseqmap}).
That $\im j^*$ is also generated by primitives of $K(G;\Q)$
follows from the 
fact that $\K(G;\Q) \isoto \H G$
is an isomorphism \emph{of Hopf algebras}%
~\cite[pf., Cor.~II.2.3]{hodgkin1967lie}.
\epf

\begin{remark}
	Much of \Cref{thm:KGKstructure} was already known;
Snaith mentioned in passing in 1971 that the K{\"u}nneth spectral sequence
beginning at $E_2 = \Tor_{R(G)}\big(\Z,R(K)\big)$,
which Hodgkin had already shown to converge to $\K(G/K)$, 
collapses at $E_2$ modulo torsion~\cite[Thm.~8.1(ii)]{hodgkin1975kunneth}\cite[p.~562]{snaith1971homogeneous}.
This implies our theorem on the level of $R(G;\Q)$-modules.
In fact, the sequence collapses even with torsion if $\pi_1 G$ is torsion-free~\cite[Thm.~2.1]{minami1975symmetric}.
With the additional assumption $\pi_1 G$ is torsion-free,
this yields explicit expressions with $\Z$ coefficients
for the nice cases enumerated in Example \ref{eg:eqf}: 
Minami proved the expected result in the case 
$R(G) \lt R(K)$ is surjective~\cite[Prop. 4.1]{minami1975symmetric}
and also, building on work of Harris,
in the case $G/K$ is a symmetric space~\cite{harris1968homogeneous,minami1975symmetric},
and of Pittie in the case where the ranks of $G$ and $K$ are equal~\cite{pittie1972homogeneous}.
\end{remark}

\section{Isotropy formality through the lens of
\textit{K}-theory}\label{sec:regularity}
This section is devoted to studying 
\isotfity of homogeneous spaces by means of \Kthy. 
Using 
the equivalence of cohomological equivariant formality and 
\Ktic \eqfity asserted by \Cref{thm:equiv} 
and the description of the \Kthy of homogeneous spaces given by 
\Cref{thm:KGKstructure}, we will give alternative proofs of 
the
\isotfity of some previously known examples. 
Furthermore, by unraveling  \Cref{thm:equiv} 
in the context of homogeneous spaces, 
we give another characterization of \isotfity in terms of representation theory. 
Before doing so, we give a $K$-theoretic version 
of the Shiga--Takahashi conditions for \isotfity 
analogous to \Cref{thm:ST++}. 
One might be tempted to simply replace the condition that 
$\rho^*\:\HG\lt \HSN$ be surjective by the condition that 
the restriction map 
$i^*\: R(G; \Q)\to R(S; \Q)^N$ be surjective,
but it turns out that this naive translation is incorrect
(see Example \ref{ex:SU4eqf}). 
Nevertheless, a slight modification admits a $K$-theoretic analogue.

\begin{proposition}\label{thm:surjectivitymodification}
	Let $(G, K)$ be a compact, connected pair and $N = \pi_0 N_G(K)$.
	\begin{enumerate}
		\item The restriction map $\rho^*: \HG \lt \HKN$ is surjective if and only if $({\HK}\quot{\HG})^N\cong\Q$.
		\item $(G, K)$ is \isotf if and only if it is formal 
		and $\big(\mspace{-1mu}{R(K; \Q)}\quot{R(G; \Q)}\mspace{-1mu}\big)^N\cong\Q$. 
	\end{enumerate}
\end{proposition}
\begin{proof}
	\begin{enumerate}
		\item \nd By the graded Nakayama lemma~\cite[Prop.~A.1.1]{neuselsmith}, 
		$\rho^*$ is surjective if and only if ${\HKN}\quot{\HG}\cong\Q$.
		But ${\HKN}\quot{\HG} = ({\HK}\quot{\HG})^N$, 
		for if $a \in \HK$ 
		represents an $N$-invariant element of $\HK \quot \HG$, 
		then so also does the $N$-average $\frac 1{|N|}\sum_{w \in N} w^* a$.
%
		\item This follows from \Cref{thm:ST++} 
			and the isomorphism $\HK \quot \HG \iso R(K;\Q) \quot R(G;\Q)$ 
			in the proof of \Cref{thm:KGKstructure}.
			\qedhere
	\end{enumerate}
\end{proof}
\subsection{Isotropy-formality of pairs arising from Lie group automorphisms}\label{isoformalgensym}
Throughout this section, 
we assume that both $G$ and $K$ are compact, connected Lie groups. 
We want to find a $K$-theoretic necessary condition for isotropy-formality of a pair $(G,K)$.
By \Cref{thm:equiv}, this is equivalent 
to the forgetful map $\K_K(G/K;\Q) \lt \K(G/K;\Q)$ being surjective.
We know from \Cref{thm:isotformal} that if 
$(G,K)$ is to be \isotf,
$G/K$ must be formal,
and hence by \Cref{thm:KGKstructure} that 
its rational \Kthy is the tensor product of a complete intersection ring
${R(K; \Q)}\quot{R(G; \Q)}$ 
and an exterior factor $\ext \wP$. 
The factor ${R(K; \Q)}\quot{R(G; \Q)}$ 
is in the image of $\a\: R(K;\Q) \to \K(G/K;\Q)$, 
and so in the $\Q$-linear span of the
classes represented by associated bundles
$G\times_K V_\rho \to G/K$ of $K$-representations $\rho$.
These bundles are equivariant with respect 
to the action induced by the left multiplication 
of $K$ on $G$,
so all elements of the factor ${R(K; \Q)}\quot{R(G; \Q)}$
admit equivariant lifts in $K_K^*(G/K; \Q)$.
To determine if the forgetful map is surjective, then, 
it remains only to see if elements of the Samelson space $\wP$ 
generating the other tensor factor
admit equivariant lifts as well. 

We digress for a while to state foundational results on the
(equivariant) \Kthy of $G$.

\begin{definition}\label{def:Grothdelta}
	\begin{enumerate}
		\item(\textup{\cite[p. 172]{hartshorne}\cite[\SS2.3]{brylinskizhang}})%
		\label{def:Grothendieck}
		Let $A$ be a ring and $B$ an $A$-algebra. 
		The module of \defd{K{\"a}hler differentials} 
		of $B$ over $A$ 
		is the quotient \defm{$\Omega_{B/A}^1$} 
		of the free $B$-module on the symbols 
		$\{db : b \in B\}$ by the $B$-submodule generated by 
		the sets
		$\{da : a \in A\}$, 
		\mbox{$\{d(b+b') - d b-d b' : b, b' \in B\}$}, and 
		$\{d(bb') - b \,db' - b' \,db : b, b' \in B\}$. 
		The commutative graded algebra of \defd{Grothendieck differentials} \defm{$\Omega^*_{B/A}$} is defined to be the exterior $B$-algebra 
		$\smash{\ext_B \Omega_{B/A}^1 =
		\bigoplus_{p=0}^\infty \bigwedge_B^p \Omega^1_{B/A}}$.
	
	\smallskip
	
	\item(\textup{\cite[{\SS}I.4]{hodgkin1967lie}\cite[\SS3]{brylinskizhang}\cite[Defs. 2.2,5]{fok2014real}})\label{def:delta} 
	Let $G$ be a compact Lie group. 
	Then $\defm{\d}\: R(G) \lt K^{-1}(G)$
	is the map which sends a complex $G$-representation $\rho$ 
	with underlying vector space $V$ to the 
	class represented by the complex of vector bundles
	\begin{align*}
		0 \to G \x \R \x V& \lt G \x \R \x V \to 0,\\
		(g, t, v)&\lmt \begin{cases}(g, t, -t\rho(g)v) \ &\text{if }t\geq 0,\\ (g, t, tv)\ &\text{if }t\leq 0.\end{cases}
	\end{align*}
	Let $G$ act on itself by conjugation and 
	$\defm{\KGad(G)}$ 
	denote the equivariant \Kthy of $G$ 
	with respect to this action.\footnote{\ 
	The superscript ``$\Ad$'' 
	emphasizes the action is by conjugation.} 
	For any subgroup $K$ of $G$ we may define a 
	map $\defm{\d_K}\: R(G) \lt \KKad(G)$ similarly: 
	$\d_K(\rho)$ is the above complex of vector bundles
	equipped with the $K$-action given by 
	$k\cdot(g, t, v) = \big(kgk\-, t, \rho(k)v\big)$.
	\end{enumerate}

\end{definition}

\begin{theorem}[Hodgkin~{\cite[p. 8 \textit{ff.}]{hodgkin1967lie}}]\label{thm:Hodgkin}
	Let $G$ be a compact, connected Lie group. 
	The map $\d\: R(G) \lt K^{-1}(G)$ is a $\Z$-linear derivation,
	or in other words satisfies
	\[
		\d(\rho_1\ox\rho_2)=\dim (\rho_1)\d(\rho_2)+\dim (\rho_2)\d(\rho_1).
	\]
	Consequently $\d(\Z) = \d \big(I(G)^2\big) = 0$
	and so $\d$ factors through a group homomorphism 
	$
		 {QR(G) \lt K^{-1}(G)}.
	$
	We have $\imd \ox \Q = P\K(G;\Q)$,
	inducing an isomorphism of Hopf algebras
	\[
		\defm{\bar\vp}\:
		\sideset{}{_{\,\Q}}\bigwedge[\imd \ox \Q] \isoto K^*(G; \Q).
	\]
\end{theorem}
	
\begin{theorem}[{\cite[\SS3]{brylinskizhang}}]\label{thm:BZstructure}
	\begin{enumerate}
		\item 
		The map $\delta_G$ in Definition \ref{def:Grothdelta}.\ref{def:delta}
		is a derivation of $R(G)$ taking values in the $R(G)$-module $K^{-1}(\Gad)$; i.e., $\delta_G$ satisfies
		\[\delta_G(\rho_1\otimes\rho_2)=\rho_1\delta_G(\rho_2)+\rho_2\delta_G(\rho_1).\]
		\item Let $G$ be a compact, connected Lie group with torsion-free fundamental group and $\Omega^*_{R(G)/\mathbb{Z}}$ the ring of Grothendieck differentials of $R(G)$ over $\mathbb{Z}$. 
		There is an $R(G)$-algebra isomorphism
	\[
	\defm\varphi: \Omega^*_{R(G)/\mathbb{Z}}\lt K_{G^{\mathrm{Ad}}}^*(G)
	\]
	defined 
	by $\varphi(\rho_V):=[G\times V]\in K_G^0(G)$
	and $\varphi(d\rho_V):=\delta_G(\rho_V)$,
	where $G$ acts on $G\times V$ by $g_0\cdot(g_1, v)=(g_0g_1g_0^{-1}, \rho_V(g_0)v)$.
	\end{enumerate}
\end{theorem}
\brmk\label{rmk:BZreduction}
We will use later the observation that 
Hodgkin's isomorphism $\bar\vp$ in \Cref{thm:Hodgkin}
is precisely the reduction mod $I(G)$
of the Brylinski--Zhang isomorphism $\vp$
in \Cref{thm:BZstructure}.
Note that for each closed subgroup $K$ of $G$, 
the map $\defm{\delta_K}\: R(G)\lt K_{K^\mathrm{Ad}}^{-1}(G)$ 
satisfies $\delta_K(\rho_1\otimes\rho_2) =
i^*(\rho_1)\delta_K(\rho_2) + i^*(\rho_2)\delta_K(\rho_1)$
and reduction mod $I(K)$ sends the image of $\d_K$
to that of $\d$.
\ermk
\begin{definition}\label{Samelsonsubspace}
	An element $\rho \in \ker \mspace{-2mu} \big(i^*\: R(G)\to R(K)\big)$
	may be written as the formal difference $\rho_1 - \rho_2\in R(G)$ 
	of complex $G$-representations whose restrictions to $K$ agree.
	If $V$ is the vector space underlying the restricted representation,
	then
	we define maps 
		$\defm{\d^{G/K}}\: \ker i^* \to K^{-1}(G/K)$ and 
		$\smash{\defm{\d^{G/K}_K}\: \ker i^*\to K^{-1}_K(G/K)}$
	sending $\rho$ to the class represented by 
	the complex of vector bundles
	\begin{align*}
		0
			\to 
		G/K\times\mathbb{R}\times V
			&\longrightarrow 
		G/K\times\mathbb{R}\times V
			\to 
		0,\\
		(gK, t, v)&\lmt \begin{cases}\big(gK, t, t\rho_1(g)\rho_2(g^{-1})v\big)\ 
										&\text{if }t\leq 0,\\ 
										(gK, t, tv)\ 
										&\text{if }t\geq 0,
						\end{cases}
	\end{align*}
	the $K$-action for the equivariant case being given by 
	$k\cdot(gK, t, v) = \big(kgK, t, \rho_1(k)v\big)$. 
\end{definition}

The signficance of this definition is that 
elements $\d^{G/K}(\rho) \in K^{-1}(G/K)$ 
admit the equivariant lifts $\d^{G/K}_K(\rho) \in K^{-1}(G/K)$
by construction
and pull back along $j\: G \lt G/K$ 
to elements $\d(\rho) \in P\K G$
by \Cref{thm:Hodgkin},
so that $\im (j^* \o \d^{G/K}) \leq \wP$ is a \Ktic analogue to 
the space $\wPz$ of Definition \ref{def:wPz}.

	\begin{notation}
		In the rest of this section, we  
		extend $\defm{i^*}$, $\defm{\d}$, $\defm{\d^{G/K}}$, and $\defm{\d_K}$
		to \textbf{rational} coefficients without further comment.
	\end{notation}


To demonstrate the utility of Theorem \ref{thm:equiv} 
we give a new proof of the (isotropy-)formality of 
compact, connected pairs arising from Lie group automorphisms.
 
\begin{theorem}\label{thm:gensympaireqf}
	Let $(G,K)$ be a \ccpair. 
	If there exists a Lie group automorphism on $G$ 
	such that the Lie algebra of the fixed point subgroup 
	coincides with the Lie algebra $\mathfrak{k}$ of $K$, 
	then
	\begin{enumerate}
		\item(\textnormal{\textup{\cite[\SS4]{terzic2001formal}%
		\cite[Prop.~4.1]{stepien2002formal}}}) 
		$G/K$ is formal and 
		\item\textnormal{(\cite{goertschesnoshari2016})} 
		$(G, K)$ is \isotf.
	\end{enumerate}
\end{theorem}
Goertsches--Noshari's proof of \isotfity consists of a series of reductions. 
First, one reduces to the case where $G$ is a simple Lie group. 
Then it is known that all automorphisms of $G$ are conjugate 
through inner automorphisms to a composition $\tau \o c_h$, 
where $\tau$ is induced by a Dynkin diagram automorphism 
with respect to a maximal torus $T$ of $G$ 
and $c_h$ is conjugation by some element $h \in T \cap  G^{\ang\tau}$%
~\cite[Lem.~5.3]{wolf1968homogeneous1}. 
If we write $K = \big(G^{\ang \s}\big)_0$ and $H = \big(G^{\ang \tau}\big)_0$, 
then within $T$, the automorphisms $\s$ and $\tau$ 
both fix exactly $S = T \cap H$. 
It turns out that $T$ is the centralizer of $S$, 
so $S$ is a maximal torus of both $K$ and $H$%
~\cite[Lem.~X.5.3, p.~492]{helgason1979differential},
and then by 
\Cref{thm:eqftorusreplacement} 
one can replace $(G,K)$ by $(G,H)$ and $\s$ by $\tau$. 
Such $\tau$ have long been enumerated, 
and the remainder of the proof involves checking equality 
of the total Betti numbers of $G/H$ and $(G/H)^S$ 
in these cases.

As for formality, Terzi\'{c} calculated $\H(G/K;\R)$ 
on a case-by-case basis and picked up formality as a corollary. 
St{\k{e}}pie\'{n} replaced $\s$ with $\tau$ as above, noted that 
symmetric spaces are already formal by work of Sullivan,\footnote{\ 
The proof is immediate from the result of \'{E}lie Cartan
that harmonic forms on a symmetric space $M$ are closed
and form a subalgebra of the de Rham algebra 
$\Omega^*(M)$ isomorphic to $\H(M;\R)$.
}
and observed the remaining case, 
$\Spin(8)/G_2 \homeo S^7 \x S^7$, 
is clearly formal.

\begin{proposition}\label{thm:gensympair}
	If $(G,K)$ is a \ccpair with $G$ simple 
	and there exists a Lie group automorphism $\sigma$ on $G$ induced by a Dynkin diagram automorphism such that the Lie algebra of the fixed point subgroup coincides with the Lie algebra $\mathfrak{k}$ of $K$, 
	then $(G, K)$ is an \isotf pair.
\end{proposition}
\begin{proof}
	We first consider the case where $G$ is simply-connected, 
	so that $R(G)$ is a polynomial ring. 
	The finite-order automorphism $\sigma$ of $G$ is induced by a graph automorphism 
	of its Dynkin diagram and the quotient graph is the Dynkin diagram of $K$. 
	Moreover, the $\rk G$ fundamental representations of $G$ 
	may be identified with the vertices of the Dynkin diagram of $G$ 
	in such a way that each element of a given $\langle \sigma\rangle$-orbit
	restricts to the same representation of $K$. 
	Let $\coprod_{k=1}^{\rk K} \{\rho_k, \s\rho_k, \cdots,\s^{j_k}\rho_k\}$
	be the partition of these fundamental representations 
	into $\langle \sigma\rangle$-orbits. 
	Then $\ker i^*$ contains
	$\rk G-\rk K$ linearly independent elements $\rho_k-\s^j\rho_k$ 
		($1\leq k\leq \rk K,\ 1 \leq j \leq j_k$) 
	The $\Q$-span of their images under $\d^{G/K}$ 
	is then a $(\rk G-\rk K)$-dimensional subspace of the Samelson space $\wP$, 
	which is itself $(\rk G-\rk K)$-dimensional 
	by the dimension inequality in \Cref{thm:purestructure}.(v), 
	so the span of these elements must be all of $\wP$. 
	It follows that $(G,K)$ is formal by \Cref{thm:purestructure}.(v). 
	Note that by construction all elements in $\wP$ admit equivariant lifts. Thus $(G, K)$ is \isotf by \Cref{thm:equiv} 
	and the discussion at the beginning of Section \ref{isoformalgensym}. 
	In the case where $G$ is not assumed simply-connected, 
	we may use \Cref{thm:eqftorusreplacement} 
	to reduce to the case where it is.
\end{proof}

Note that the above proof does not use the classification 
of generalized symmetric spaces. 
Moreover it shows formality by constructing elements of
the Samelson space without directly invoking \Cref{thm:isotformal}.

\begin{remark}\label{rmk:betahistory}
	\begin{enumerate}
		\item  Recall that the odd \Kthy functor $K^{-1}$ 
		is represented by the infinite unitary group 
		$\U(\infty)\ceq \colim \U(n)$. 
		Hodgkin~\cite[pf., Cor.~II.2.3]{hodgkin1967lie} 
		considered the map $\b: R(G)\ \lt K^{-1}(G)$ 
		sending a complex $G$-representation $\rho$ 
		to the homotopy class of 
		$\smash{G \os{\rho}{\to}\U(n)\inc \U(\infty)}$.
		This is in fact the same as $\d$,
		but we use $\d$ because
		the equivariant lift $\d_G$ is slightly easier to describe than 
		the corresponding $\b_G$.
		The construction $\d^{G/K}$ appears
		in Hodgkin's work~\cite[\SS10]{hodgkin1975kunneth} 
		on his equivariant K{\"u}nneth spectral sequence.
				
		\smallskip
				
		\item The middle map of the complex of vector bundles
		in Definition \ref{def:delta}
		is a correction of the original definition due to Brylinski and Zhang,
		which Brylinski later noted was flawed and corrected.
		The result itself is unaffected.
		Discussion of the error and correction 
		can be found in work of the second-named author~\cite[Rmk.~2.6]{fok2014real}.
	\end{enumerate}
\end{remark}

\subsection{A representation-theoretic characterization of \isotfity}
In this section we assume that $G$ is a compact, connected Lie group with torsion-free fundamental group unless otherwise specified.


%

	\begin{notation}
		In the rest of this section, we write 
		$\defm{R}$ for the image of the restriction map 
		$i^*\: R(G;\Q) \lt R(K;\Q)$
		and $\defm I$ for its augmentation ideal $i^* I(G;\Q)$,
		and implicitly extend $\d$ and $\d_K$ to $\Q$ coefficients.
		We denote by $\defm{\mathcal{M}}$
		the $R(K;\Q)$-submodule $\linspan_{R(K;\Q)}\,  \d_K(\ker i^*)$
		of $K^{-1}_{\Kad}(G;\Q)$.
	\end{notation}
	
	\begin{lemma}\label{repring}
		Let $(G,K)$ be a \ccpair and $R = \im i^*$ as above. Then
		\begin{enumerate}	
			\item\label{fgintrepring} both $R(G; \Q)$ and $R$ are 
			integral domains finitely generated as $\Q$-algebras,
			\item\label{rankofrepring} $\Kdim R(G)=\rk G$, and
			\item\label{freeKahlerRG} if $\pi_1 G$ is torsion-free,
				$\smash{\Omega^1_{R(G;\Q)/\Q}}$
			is a free $R(G;\Q)$ module of rank equal to $\rk G$.
		\end{enumerate}
	\end{lemma}
	\begin{proof}
		Let $T$ be a maximal torus of $G$.
		Then restriction of representations induces an isomorphism
		$R(G) \isoto R(T)^W$~\cite[{\SS}4.4]{atiyahhirzebruch}.
		\benum
		\item
		As $R(T;\Q)$ is a Laurent polynomial ring on $\rk G$ 
		generators over $\Q$, it follows 
		$R(G;\Q) \leq R(T;\Q)$
		and $R \leq R(K;\Q)$
		are integral domains.
		As $R(T;\Q)$ is finitely generated over $\Q$,
		so are $R(G;\Q)\iso R(T;\Q)^W$~\cite[Thm.~2.1.4]{neuselsmith} 
		and its quotient $R$. 
		In fact, the connectedness hypothesis is not needed~\cite[Cor.~3.3]{segal1968representation}.
		\item 
		As $R(T; \Q)$ 
		is integral over $R(G; \Q) \iso R(T; \Q)^W$~\cite[Ex. 5.12]{atiyahmacdonald}, 
		$\Kdim R(G; \Q) = \Kdim R(T; \Q)$ 
		by going-up and lying-over. 
		But it is clear $\Kdim R(T; \Q)=\rk G$.\footnote{\ 
		Explicitly, $R(T;\Q) = \Q[t_1^{\pm 1},\ldots,t_{\rk G}^{\pm 1}]$ 
		is a localization of $\Q[t_1,\ldots,t_{\rk G}]$ 
		and hence $\Kdim R(T;\Q) \leq \Kdim \Q[t_1,\ldots,t_{\rk G}] = \rk G$
		on the one hand~\cite[Ex., p.~121]{atiyahmacdonald},
		and on the other 
		$(0) < (t_1 - 1) < \cdots < (t_1-1,\ldots,t_{\rk G} -1)$
		is a chain of prime ideals of length $\rk G$ in $R(T;\Q)$.
	}
		\item
		 When $\pi_1 G$ is torsion-free, $R(G; \Q)$ 
		 is the tensor product of a polynomial ring %
		 and a Laurent polynomial ring
		 ~\cite[Prop.~11.1]{hodgkin1975kunneth},
		 and particularly, is a localization of a polynomial ring.
		 But localization commutes with 
		 K{\"a}hler differentials~\cite[Cor.~12.2.16]{CRing},
		 and if $k$ is a field,
		 then $\smash{\Omega^1_{k[x_1,\ldots,x_n]/k}}$
		 is the free $k[x_1,\ldots,x_n]$-module on $dx_1,\ldots,dx_n$%
		~\cite[12.2.8]{CRing}.
		\qedhere
		\eenum
	\end{proof}


\begin{lemma}\label{thm:affinesmooth}
	Let $k$ be a field of characteristic zero, 
	$A$ an integral domain finitely generated as a $k$-algebra,
	and $L$ the field of fractions of $A$.
	Then $\rk_A \Omega^1_{A/k} \ceq \dim_L L \ox_A \Omega^1_{A/k} = \Kdim A$.
\end{lemma}
\begin{proof}
By Noether normalization, 
$A$ is integral over some polynomial $k$-subalgebra $B$,
which by going-up and lying-over
has equal Krull dimension.
Now $\Kdim B$ is the transcendence degree of $L$ over $k$
since $B$ is a polynomial ring over $k$ and 
its field of fractions is $L$.
But $\mathrm{tr\,deg}_{k} L = \dim_L \Omega^1_{L/k}$~\cite[Cor.~12.3.5]{CRing},
and $\Omega^1_{L/k} \iso L \ox_A \Omega^1_{A/k}$ 
since localization commutes with 
taking K{\"a}hler differentials~\cite[Cor.~12.2.16]{CRing}.
\end{proof}

	\begin{proposition}\label{thm:descriptSamelson}
		The module $\mathcal{M}$ 
		is of rank equal to $\rk G - \rk K$. 
	\end{proposition}
	\begin{proof}
	    We have the 
	    standard exact sequence~\cite[Prop. 8.4A, Chap. II]{hartshorne}
	    of $R$-modules
		\begin{eqnarray}\label{eq:secondexactseq}
			\quotientmed{\keri }{(\keri)^2}
				\,\lt\,
			\xt{R(G;\Q)}R{\Omega^1_{R(G; \Q)/\Q}}
				\,\lt\,
			\Omega^1_{R/\Q}
				\,\to\,
			0,
		\end{eqnarray}
		where the first map sends $\rho + (\keri)^2 \lmt 1 \ox d\rho$.
		Extending coefficients gives an exact sequence
		\[
			\xt {R(G;\Q)}{R(K;\Q)}{\frac{\keri }{(\keri)^2}}
				\,\lt\, 
			\underbrace{\xt {R(G;\Q)}{R(K;\Q)}{\Omega^1_{R(G; \Q)/\Q}}}_{\defm{\mathcal N}}
				\,\lt\,
			\underbrace{\xt R{R(K;\Q)}{\,\Omega^1_{R/\Q}}}_{\defm{\mathcal P}}
				\,\to\,
			0.
		\]
		Thus the kernel of $\cal N \lt \cal P$ is $\mathcal M$.
		By \Cref{repring}.\ref{freeKahlerRG},
		we have $\rk_{R(K;\Q)} \mathcal N = \rk G$
		and as for $\mathcal P$, we find 
		\[
 			   \rk_{R(K; \Q)}\mathcal P
			= \rk_{R}\Omega^1_{R/\Q}
			= \Kdim R 
			= \Kdim R(K;\Q)
			= \rk K,
		\]
where the second equality is
\Cref{thm:affinesmooth},
the third follows from going-up
since $R(K;\Q)$ is integral over $K$~\cite[Prop.~3.2]{segal1968representation}%
\cite[Rmk., p.~60]{atiyahmacdonald},
and the last is \Cref{repring}.\ref{rankofrepring}.
By exactness, 
$\rk_{R(K;\Q)} \mathcal M = \rk_{R(K;\Q)} \mathcal N - \rk_{R(K;\Q)} \mathcal P$.
\end{proof}
\bprop\label{thm:KSamelsondeltaformal}
	Let $(G,K)$ be a \ccpair.
	Then $(G,K)$ is \isotf if and only if 
	$\dim q (\ker i^*) = 
	\dim \delta(\ker i^*) =
	\dim (\im \d^{G/K}) = 
	\rk G - \rk K$.
\eprop
\bpf
	From (\ref{eq:cohomseqmap}), the map $\text{ch}_G$ restricts to 
	the map $\keri \lt \ker \rho^*$.
	Since $R(G;\Q)\compl \iso \HG$ by \Cref{thm:eqcherniso},
	applying \Cref{thm:tensorreduction}
	in the case $A = R(G;\Q)$, $M = I(G;\Q)$, and $N = \Q$
	shows that $QR(G;\Q) \lt Q\HG$ is also an isomorphism
	restricting to an injection $q \keri \lt  q \ker \rho^*$
	of indecomposable images.
	The result will follow from
	\Cref{thm:HSamelsondeltaformal}, \Cref{thm:equiv}, 
	and the discussion in the first paragraph of this subsection if we can show 
	this last map is also surjective. 
	For this, recall the $I(K)$-adic and $I(G)$-adic topologies 
	on $R(K)$ agree~\cite[Cor.~3.9]{segal1968representation}
	and complete the short exact sequence
	$
		0 \to \ker i^* \to R(G;\Q) \to R(K;\Q) \to 0 
	$ 
	at $\fa = I(G;\Q) \idealneq R(G;\Q) = A$ 
	to conclude $\ker \rho^* = (\keri)\compl$.
	Writing $\fb = \keri \leq \fa$,
	our task is to see
	$(\fb + \fa^2)/\fa^2 \lt (\wh\fb + \wh\fa^2)/\wh\fa^2$
	is surjective.
	It is equivalent to see 
	$
	\fb 
		 \lt
	(\wh\fb + \wh\fa^2)/\wh\fa^2 
	$
	is surjective,
	but that map factors as
	\[
	\fb 
		\longepi 
	\fb/\fa \fb 
		\isoto 
	\wh\fb/\wh\fa\wh\fb 
		\longepi 
	\quotientmed{\wh\fb\,}{\,\wh\fb \inter \wh\fa^2}
		\isoto 
	\quotientmed{\wh\fb + \wh\fa^2\,}{\,\wh\fa^2}.	
	\]
	Each factor is obviously surjective except, perhaps, the second,
	which follows from the string of standard isomorphisms
	\[
	\fb/\fa\fb \iso 
	\xt{A} \fb {A /\fa} \iso
	\xt{A} \fb {\wh A / \wh\fa}\iso
	\xt{\wh A}{\wh\fb}{\wh A / \wh\fa} \iso
	\wh\fb/\wh\fa\wh\fb;
	\]
	the penultimate isomorphism comes from \Cref{thm:tensorreduction}.
\epf
\bpf[Proof of \Cref{thm:regularity}]
	We start by rephrasing regularity of $R$ at $I$ 
	in numerical fashion.
	By \Cref{thm:eqftorusreplacement} and the assumption made at the beginning of this subsection, 
	we may assume that $G=\tG$. 
	By \Cref{thm:affinesmooth},
	we know $\dim_{R_{(0)}} {R_{(0)}} \ox _R{\,\Omega^1_{R/\Q}}
	= \Kdim R$,
	and 
	since $R$ is an integral domain finitely generated over $\Q$
	by Lemma \ref{repring}.\ref{fgintrepring}
	and $I$ is a maximal ideal, 
	this is also $\Kdim R_I$~\cite[Cor.~11.27]{atiyahmacdonald}. 
	Now, $R$ is regular at $I$ if and only if~\cite[Cor.~3.7]{CRing} 
	these numbers equal
	 \[
		 \dim_\Q \xt R {(R/I)} {\,\Omega^1_{R/\Q}}
			 = 
		 \dim_\Q \xt R {\Q\,}{\,\Omega^1_{R/\Q}}.
	\]
		Tensoring the exact sequence (\ref{eq:secondexactseq}) over $R$ 
		with $R_{(0)}$ and $\Q$, we find respectively
		\quation{\label{eq:dimeqprelim}
			\xt R{R_{(0)}}{\,\Omega^1_{R/\Q}} 
				\iso 
			\frac{ {R_{(0)}} \ox_ {R(G;\Q)}{\,\Omega^1_{R(G;\Q)/\Q}}}
				{R_{(0)} \ox_{R(G;\Q)} d(\ker i^*)}
			\quad\ \
				\mbox{and}
			\quad\ \ 
			\xt R {\Q\,}{\,\Omega^1_{R/\Q}}
				\iso
			\frac{{\Q\,}\ox_ {R(G;\Q)}{\,\Omega^1_{R(G;\Q)/\Q}}}
			{\Q \ox_{R(G;\Q)} d(\ker i^*)}.
		}
		Since $\Omega^1_{R(G; \Q)/\Q}$ is a free $R(G; \Q)$-module 
		by \Cref{repring}.\ref{freeKahlerRG},
		the numerators in (\ref{eq:dimeqprelim}) 
		are of equal dimension over their respective scalar fields,
		and hence we have regularity if and only if 
		\[
			\dim_{R_{(0)}} \xt{R(G; \Q)}{R_{(0)}} {d(\ker i^*)}
			=
			\dim_\Q \xt{R(G; \Q)}{\Q} {d(\ker i^*)}.
		\]
		Now, using 
		the isomorphisms of
		Theorems \ref{thm:Hodgkin} and \ref{thm:BZstructure}
		and the observation of Remark \ref{rmk:BZreduction}
		linking them,
		the left- and right-hand sides may respectively
		be identified with base extensions of the images of
		$\d_K(\ker i^*)$ and $\d(\ker i^*)$, respectively,
		so the equation is
		\[
			\dim_{R_{(0)}}\linspan_{R_{(0)}} \delta_K(\keri)
			=
			\dim_\Q \linspan_\Q \delta(\keri).
		\]
		Since $R(K; \Q)_{(0)}$ is a field extension of $R_{(0)}$, 
		the left-hand side is $\rk_{R(K;\Q)} \mathcal M$,
		which is equal to $\rk G - \rk K$
		by \Cref{thm:descriptSamelson},
		so finally $R$ is regular at $I$
		if and only if  $\dim_\Q \linspan_\Q \d(\keri) = \rk G - \rk K$.
		But by \Cref{thm:KSamelsondeltaformal},
		this happens if and only if $(G,K)$ is \isotf.
\epf
		
%
%

\begin{remark}\label{rmk:onefellswoop}
	There is another way of interpreting the condition that $i^*R(G; \Q)$ be regular at $I$. 
	As mentioned in the proof of Lemma \ref{repring}.\ref{freeKahlerRG}, $R(G; \Q)$ is the ring
	$\smash{\Q[\rho_1, \cdots, \rho_i, t_1^{\pm 1}, \cdots, t_{l-i}^{\pm 1}]}$. 
	If we let $\smash{\overline{\rho}_j=\rho_j-\dim \rho_j}$
	and $\smash{\overline{t}_k^{\pm 1}=t_k^{\pm 1}-1}$ 
	be the ``reduced representations'' 
	and $\ker \wt{\imath}^*$ is minimally generated by 
	$(k_1, \cdots, k_p)$, 
	then $R = \im i^*$ is regular at $I$ if and only if each of 
	the generators
	$k_1, \cdots, k_p$, when written as a polynomial in the reduced representations, has nonzero linear terms, i.e., is not in $I^2$. 
\end{remark}

\brmk\label{rmk:regularH}
Write $\defm{\wh R}$
for the completion  of $R$ at $I$ 
and $\defm{\wh I}$ for its augmentation ideal.
Then the regularity of $R$ at $I$
is equivalent to the regularity of 
$\wh R$ at $\wh I$ since $\wh{R}_{\wh I} \iso \wh{R_I}$
and \cite[Prop.~11.24]{atiyahmacdonald}
a Noetherian local ring is regular if and only if its completion is regular.
From \Cref{thm:eqcherniso},
since the $I(K)$-adic and $I(G)$-adic topologies on $R(K)$
agree~\cite[Cor.~3.9]{segal1968representation}, 
we see ${\im(H_G^{**} \lt H_K^{**})}\iso \wh R$,
so the regularity condition in \Cref{thm:regularity}
can also be phrased in terms of cohomology.
\ermk

	\Cref{thm:regularity} allows us to give a uniform proof of \isotfity (and hence formality by \Cref{thm:isotformal}) of some classes of homogeneous spaces in Examples \ref{eg:eqf}. 
\begin{proof}[Proofs for Example \ref{eg:eqf}]
\begin{enumerate}
	\item If $(G, K)$ is an equal-rank pair of compact, connected Lie groups, 
	so is $(\tG, \tK)$. 
	The restriction map $\wt{\imath}: R(\tG; \Q)\to R(\tK; \Q)$ is injective, 
	so the image $R \iso R(\tG; \Q)$ is a polynomial ring 
	tensored with a Laurent polynomial ring
	and thus regular at $I$.

	\item Let $(G,K)$ be a \ccpair such that $\H G \lt \H K$
	is surjective. 
	Then~\cite[Cor., p.~179]{borelthesis} $H_G^{**} \lt H_K^{**}$ is surjective
	as well, so $\wh R = H_K^{**}$, which is regular
	at its augmentation ideal $\wh I$ since it is a power series ring,
	so by Remark \ref{rmk:regularH}, we see $(G,K)$ is \isotf.
	
	\item Let $(G, K)$ be a generalized symmetric pair. 
	As in the discussion following \Cref{thm:gensympaireqf},
	we may reduce to the case $G$ is simple and simply-connected
	and the Lie group automorphism of $G$ is induced by a 
	graph automorphism of its Dynkin diagram. 
	Our proof of \Cref{thm:gensympair} 
	shows that $\ker\! \big(i^*\: R(G; \Q)\to R(K; \Q)\big)$ 
	is an ideal generated by linear combinations of reduced fundamental representations of $G$, 
	and by Remark \ref{rmk:onefellswoop}, 
	$\im i^*$ is regular at $I$. 
	\qedhere
\end{enumerate}
\end{proof}

\subsection{Some examples}\label{sec:endeg}
In the following examples, $G$ is a special unitary group, 
and the torus subgroup $S$ we consider is one-dimensional.
Thus the fundamental group of $G$ is trivial 
and we do not need to consider a central cover,
and $R(S)$ is a principal ideal domain, 
so that $R(S) \quot R(G)$ is a complete intersection ring
and $G/S$ is formal
by Theorems \ref{thm:KGKstructure} and \ref{thm:purestructure}.

\begin{example}\label{ex:SU4eqf}
Let $G=\SU(4)$ and $S = \big\{\mn\diag(z, z^{-1}, z^2, z^{-2}) : z \in S^1\big\}$.
By Example \ref{eg:eqf}.\ref{eg:circle}, 
since $S$ is reflected, the pair $(G,S)$ is \isotf. 
We can also show \isotfity 
by verifying the Shiga--Takahashi criteria of \Cref{thm:ST++}:
the restriction map $\Q[c_2,c_3,c_4] \iso \HG \to \HS \iso\Q[s]$ 
assigns the universal Chern classes as elementary symmetric polynomials 
in $(s,-s,2s,-2s)$:
\[
		c_2 \lmt -5s^2, \qquad\qquad
		c_3 \lmt 0, \qquad\qquad
		c_4 \lmt 4s^4,
\]
and so the image is $\Q[s^2] \iso \HSN$ 
and the cokernel $\HS\quot\HG \iso \Q[s]/(s^2)$.
	
	We could also use \Cref{thm:regularity} to show \isotfity of $(G,S)$. 
	If $\s_4$ is the defining representation of $\SU(4)$,
	the reduced fundamental representations 
	$x\ceq \sigma_4-4$, \ 
	$y\ceq \varbigwedge^2\sigma_4-6$, and 
	$z\ceq \varbigwedge^3\sigma_4-4$
	generate the augmentation ideal $I(G)$ 
	of $R(G) \iso\Z[x, y, z]$.
	If we write $R(S) = \Z[t,t\-]$ and $a = t + t\- -2$,
	then $R(S)^N = \Z[a]$.
	Under the restriction $i^*\: R(G) \lt R(S)^N$ of representations,
	the generators map as
	\[
		x \lmt a^2 + 5a,\qquad\qquad
		y \lmt a^3 + 6a^2 + 10a,\qquad\qquad
		z \lmt a^2 + 5a,
	\]
	so the image $\imi$ is  
		$
			\Z[a^2 + 5a,a^3 + 6a^2 + 10a]
		$.
	This is a proper subring of $\Z[a] = R(S)^N$, 
	showing the naive modification of the 
	Shiga--Takahashi criterion of \Cref{thm:ST++} does not hold.
	However, the modification in \Cref{thm:surjectivitymodification}
	does, for
%
	as predicted by \Cref{thm:KGKstructure}, 
	we have an isomorphism
	\begin{align*}
		R(S;\Q) \quot R(G;\Q)
			\,&\, =
		\Q[t,t\-] /(a)
			\\ \,&\, =
		\quotientmed{ 	\Q[t,a]	\,	}{	\big(a, (t-1)^2 - at	\big)	}
			\\ \,&\, \iso
		\Q[t] / (t-1)^2	
			\,\iso\, 	
		\HS \quot \HG
	\end{align*}
	and $(\HS \quot \HG)^N \iso \big(\Q[s]/(s^2)\big)^N = \Q$.
Computing the resultant of the equations
\[
a^2+5a-x = 0
\qquad\mbox{and}\qquad
a^3 + 6a^2 + 10a - y = 0
\]
with respect to $a$,
we find the kernel of $i^*$ is 
$(x-z, -x^3-14x^2+3xy-50x+y^2+25y)$,
so $R \iso \Q[x, y, z] / \keri $
is regular at $I$ by Remark \ref{rmk:onefellswoop}, and by Theorem \ref{thm:regularity} $(G, S)$ is isotropy-formal. 

In fact, one can easily check that the $\d^{G/S}$-images of $x - z$ and $-x^3-14x^2+3xy-50x+y^2+25y$ 
are linearly independent, so \Cref{thm:KSamelsondeltaformal}
again shows $(G,S)$ is isotropy-formal and hence formal by Theorem \ref{thm:isotformal}. We have that
\[
	\K(G/S;\Q) 
	\,\iso\,
	\frac{\Q[t-1]\,}{\,(t-1)^2} 
	\,\tensor
		\sideset{}{_{\,\Q}\!}\bigwedge\big[\d^{G/S}(x-z), \d^{G/S}(-x^3-14x^2+3xy-50x+y^2+25y)\big].
\]
The displayed exterior generators of course admit equivariant lifts.
\end{example}
	
\begin{example}\label{ex:SU3noteqf}
	Let $G = \SU(3)$ 
	and $S = \big\{\!\diag(z, z, z^{-2}) : z \in S^1\big\}$.
	As $S$ is not reflected,
	$(G,S)$ cannot be \isotf
	by Example \ref{eg:eqf}.\ref{eg:circle}.
	Alternately, we can see this using \Cref{thm:regularity} 
	and showing $i^* R(G;\Q)$ is not regular at $I$.
	Let $\s_3$ be the defining representation of $G$ so that
	$x = \sigma_3 - 3$ and $y= {\varbigwedge}^2	\sigma_3 - 3$ generate $I(G)$
	and $R(G) = \Z[x,y]$.
	Then $i^*$ takes 
	\[
		x \lmt	2t+t^{-2} -3,\qquad\qquad
		y \lmt 2t\-+t^2-3.
	\]
	Computing the resultant,
	one finds
	$\ker i^* = ( 4x^3+4y^3-x^2y^2-6x^2y-6xy^2+27x^2+27y^2-54xy)$, 
	which nontrivially intersects $I(G)^2 = (x,y)^2$,
	so $R$ is not regular at $I$.
\end{example}

{\footnotesize\bibliography{bibshort.bib} }

\smallskip

\nd\footnotesize{\textsc{ }\\
	\url{jeffrey.carlson@tufts.edu}
}

\medskip

\nd\footnotesize{\textsc{Xi'an Jiaotong-Liverpool University,\\
111 Ren’ai Road, Suzhou Industrial
Park,\\Suzhou, Jiangsu Province 215123, China}\\
	\url{chikwong.fok@xjtlu.edu.cn}
}

\end{document}